\documentclass[12pt]{article}
\usepackage{amssymb,amsthm,amscd,amsmath,graphicx,color}

\newcommand{\R}{{\mathbb R}}

\newcommand{\Z}{{\mathbb Z}}
\newtheorem{theorem}{Theorem}[section]
\newtheorem{corollary}[theorem]{Corollary}
\newtheorem{lemma}[theorem]{Lemma}
\newtheorem{proposition}[theorem]{Proposition}

\newtheorem{remark}[theorem]{Remark}

\begin{document}

\thispagestyle{empty}
\begin{flushright} \rm J. Geo. Phys. {\bf
  56/9} (2006), 1875-1892 \end{flushright}
\par\bigskip\par
\vfill
\begin{center}
{\bfseries\Large
Two-dimensional Lagrangian singularities and bifurcations of gradient lines II}
\par\addvspace{20pt}
{\sc G. Marelli}
\par\medskip
Department of Mathematics, Kyoto University,\\
Kitashirakawa, Sakyo-ku, Kyoto 606-8502, Japan\end{center}
\vfill
\begin{quote} \footnotesize {\sc Abstract.}
Motivated by mirror symmetry,
we consider the Lagrangian fibration $\R^4\rightarrow\R^2$ and
Lagrangian maps $f:L\hookrightarrow \R^4\rightarrow \R^2$, 
exhibiting an unstable singularity, and study 
how 
the bifurcation locus of gradient lines, the integral curves of 
$\nabla f_x$, for $x\in B$,
where $f_x(y)=f(y)-x\cdot y$, changes when $f$ is slightly perturbed.
We consider the cases when $f$ is the germ of a fold, of a cusp
and, particularly, of an elliptic umbilic.
\end{quote}
\vfill
\leftline{\hbox to8cm{\hrulefill}}\par
{\footnotesize
\noindent Supported by a JSPS Postdoctoral fellowship
\par
\noindent\emph{2000 Mathematics Subject Classification:} 37G25, 53D12,
70K60 \par
\noindent\emph{E-Mail address:} {\tt marelli@kusm.kyoto-u.ac.jp}
}
\eject\thispagestyle{empty}\mbox{\ \ \ }\newpage\setcounter{page}{1}

\section{Introduction}
This is the second of two papers motivated by the problem
of quantum corrections in mirror symmetry. More precisely, in the first part
\cite{M}, we considered the torus fibration $T^4\rightarrow T^2$ and a 
Lagrangian map $f:L\hookrightarrow T^4\rightarrow T^2$ exhibiting some 
unstable singularity, and studied 
how this singularity break when $L$ is slightly perturbed. We defined
then the gradient lines of $f$ and examined some of their properties.
Here instead we want to analyse, in a neighbourhood of the caustic of $f$, 
how the bifurcation locus of gradient lines
changes when $L$ is perturbed.
Since the problem is local, we can simply consider the Lagrangian fibration
$\R^4\rightarrow\R^2$.
The study of the bifurcation locus of gradient lines for perturbations of $f$ 
is harder than the similar one, examined in \cite{M},
regarding the caustic: indeed, since 
the caustic is the set of critical values of $f$,
the matter is, in a sense, local, while it is global
when studying 
bifurcations of solutions of $\nabla f_x=0$, since global aspects
of the flow of a vector field are involved. 
The problem is
relatively easy for the fold and the cusp, but
it becomes much harder for perturbations of 
the elliptic umbilic, because, being its 
bifurcation locus non-empty, its complexity increases
considerably when a small perturbation is added. The main result of this paper
is theorem \ref{finalresult}, which is concerned with the bifurcation diagram
of a small perturbation of the elliptic umbilic. The
study of all these cases should help in
drawing some ideas about the mutual positions of
bifurcation lines and caustic.

\par\smallskip
{\bf Acknowledgements.} I wish to thank K.~Fukaya, whose suggestions
and 
help were decisive for
the achievement of all the results here expounded. 

I am thankful to JSPS (Japan Society for the 
Promotion of Science) which awarded me with a JSPS postdoctoral
fellowship at 
Kyoto University,
where this paper was written.

\section{The fold}
The dynamical system (10) in \cite{M},
when $f$ is the generating function (2) in \cite{M}
of the fold in dimension 2, takes the form
\begin{equation}
\label{dsf}
\left\{ \begin{aligned}
\frac{dy_1}{dt} & =y_1^2-x_1\\
\frac{dy_2}{dt} & =y_2-x_2
\end{aligned}
\right.
\end{equation}

\begin{proposition}
${\cal B}=\emptyset$.
\end{proposition}

\begin{proof}
The caustic has equation $x_1=0$.
The vector fields (\ref{dsf}) has respectively two, one or no critical points,
depending on whether $x_1>0$, $x_1=0$ or $x_1<0$.
If $x_1>0$ the critical points are
$(\sqrt{x_1},x_2)$ and $(-\sqrt{x_1},x_2)$. Linearizing
the vector field 
in a neighbourhood of these points, we find out that $(\sqrt{x_1},x_2)$
has two positive eingenvalues, so it is an unstable node, while
$(-\sqrt{x_1},x_2)$ has a positive and a negative eigenvalue, so it is
a saddle. Equations (\ref{dsf}) can be easily solved: there is a 
gradient line from the node to the saddle, which is generic and
stable, and
whose image is the 
segment with the node and the saddle as extremes. 
At $x_1=0$, on the
caustic, the node and the saddle glue together in a 
saddle-node.
\end{proof}

\section{The cusp}

For the cusp
\begin{equation}
\label{gfc}
f(y_1,y_2)=\frac{1}{4}y_1^4+y_1^2y_2+\frac{1}{2}y_2^2
\end{equation}
the dynamical system (10) in \cite{M}
has the form
\begin{equation}
\label{dscusp}
\left\{ \begin{aligned}
\frac{dy_1}{dt} & = y_1^3+2y_1y_2-x_1\\
\frac{dy_2}{dt} & = y_1^2+y_2-x_2
\end{aligned}
\right.
\end{equation}
and the caustic is the semicubical parabola 
$$|x_1|=\frac{4}{3}\sqrt{\frac{2}{3}}x_2^{3/2}$$
defined for $x_2\geq0$.

\begin{proposition}
\label{cuspbifu}
The bifurcation diagram of a small perturbation of 
the cusp, supported on a sufficiently small 
compact set $W$, has the following features,
outlined in figure \ref{bifcusp2}.1:\\ 
1. for 
$|x_1|>\frac{4}{3}\sqrt{\frac{2}{3}}x_2^{3/2}$ 
(``outside''
the caustic), $\nabla f_x$ has only one critical point: a saddle;\\
for 
$|x_1|<\frac{4}{3}\sqrt{\frac{2}{3}}x_2^{3/2}$ (``inside'' the
caustic), $\nabla f_x$ has three critical points: two
saddles, $s_1$ and $s_2$, and a node $n$;\\
for 
$x_1=\frac{4}{3}\sqrt{\frac{2}{3}}x_2^{3/2}$ (the $l_1$
branch of the caustic),
$\nabla f_x$ has two critical points: the saddle $s_2$ and the
saddle-node $ns_1$;\\ 
for 
$x_1=-\frac{4}{3}\sqrt{\frac{2}{3}}x_2^{3/2}$ (the $l_2$
branch of the caustic),
$\nabla f_x$ has two critical points: the saddle $s_1$ and the
saddle-node $ns_2$;\\ 
at $(0,0)$ (the vertex of the cusp), $\nabla f_x$ has a degenerate
critical point $ns_1s_2$;\\
2. if non-empty, ${\cal B}$ lies inside the caustic; moreover there
exist a neighbourhood $U$ of
the vertex of the cusp inside the caustic such that 
$U\cap{\cal B}=\emptyset$; in $W$, ${\cal B}_{ij}$, 
if non-empty, can contain lines,
half-lines with origin at a point of $l_i$, 
segments with both the extremes on $l_i$ or immersed $S^1$'s;
${\cal B}_{ij}\cap{\cal B}_{ji}=\emptyset$, while two components 
${\cal B}_{ij}^1$ and ${\cal B}_{ij}^2$ of
${\cal B}_{ij}$ can intersect provided the saddle-to-saddle separatrices
$\gamma_{s_is_j}$ of the vector fields corresponding to
points of ${\cal B}_{ij}^1$ and ${\cal B}_{ij}^2$ are obtained as
intersection of the same pair of 
separatrices of $s_1$ and $s_2$.
\end{proposition}

\setlength{\unitlength}{1cm}
\begin{picture}(5,5.5)
\label{bifcusp2}
\thicklines

\qbezier(1,3)(3,3)(7,5)
\qbezier(1,3)(3,3)(7,1)

\thinlines
\qbezier(3.4,3.5)(3.8,3.2)(4.2,3.2)
\qbezier(4.2,3.2)(5.5,3.2)(5.5,3.7)
\qbezier(5.5,3.7)(5.5,4.1)(6,4.1)
\qbezier(6,4.1)(6.5,4.1)(6.5,3.7)
\qbezier(6.5,3.7)(6.5,3.4)(7,3.4)

\qbezier(4.2,3.8)(4.2,3)(5.2,3)
\qbezier(5.2,3)(6.4,3.3)(6,4.5)

\qbezier(2.7,2.8)(2.7,3)(3,3)
\qbezier(3,3)(3.3,3)(3.3,2.6)

\qbezier(3.7,2.4)(3.7,3)(4.2,3)
\qbezier(4.2,3)(4.5,3)(4.5,2.7)
\qbezier(4.5,2.7)(4.5,2.4)(5,2.4)
\qbezier(5,2.4)(6,2.4)(6,2.8)
\qbezier(6,2.8)(6,3.1)(7,3.1)

\qbezier(4.5,2.1)(4.5,2.8)(5.7,3)
\qbezier(5.7,3)(6.5,3)(7,2)

\qbezier(5.3,1.8)(5.7,2.1)(6.4,2.1)
\qbezier(6.4,2.1)(6.8,2.1)(6.8,1.1)

\qbezier(3.6,2.8)(3.4,2.9)(3.6,3)
\qbezier(3.6,2.8)(3.8,2.9)(3.6,3)

\qbezier(4.6,3.4)(4.4,3.5)(4.6,3.7)
\qbezier(4.6,3.4)(4.8,3.5)(4.6,3.7)

\end{picture}
\begin{center}
$Fig.~\ref{bifcusp2}.1:~The~bifurcation~diagram~of~a~small~perturbation
~of~the~cusp$
\end{center}
Note that, while Lagrangian equivalent maps have diffeomorphic caustics,
their bifurcation loci in general are not diffeomorphic. Note also that
no bifurcation line has the vertex of the cusp as limit point.

\begin{proof}
1. All the statements follow from a direct computation
of roots of $\nabla f_x$ and from the study of the sign of 
eigenvalues of the linearization of $\nabla f_x$ in a neighbourhood of
its critical points.\\
2. Bifurcation points, if they exist, lie only inside the caustic, where
at least two saddles exist. Since, if non-empty and far from the caustic, 
the components of 
${\cal B}_{ij}$ are  
immersed submanifold of codimension 1, 
it follows that they can be either 
lines, half-lines with origin on the caustic, 
segments with extremes on the caustic, or
immersed $S^1$'s.\\ 
Suppose $\overline{{\cal B}_{ij}}\cap l_j\neq\emptyset$ 
and take $p\in\overline{{\cal B}_{ij}}\cap l_j$,
then, by proposition (4.14) in \cite{M}, 
there exists
an open subset $V$ such that $p\in\partial V$ and
the phase portrait of $\nabla
f_x$, for $x\in V$, does
not exhibit the
gradient line $\gamma_{ns_j}$; when $x$ moves along a path in $V$ ending at 
$p$, $n$ and 
$s_j$ glue together at
$p$ forming the degenerate critical point $ns_j$, and 
so in the phase portrait of 
$\nabla f_p$ there would be two gradient lines with opposite directions
joining $ns_j$ and $s_i$, which is not possible, as shown in the proof of
corollary (4.13) in \cite{M}.\\
Suppose the vertex $v$ of the cusp belongs to $\overline{{\cal B}}$, 
then moving along paths ending at $v$ contained into different 
connected components 
determined by ${\cal B}$, 
the phase portrait of $\nabla f_v$ would depend on these
paths, giving a contraddiction.\\
The statement about the intersection of components of ${\cal B}$ is a 
consequence of corollary (4.18) in \cite{M}.\\
Finally, if $W$ is sufficiently small,
theorem (4.30) in \cite{M} 
ensures in $W$ 
the existence of a generating function having
a bifurcation diagram with such features.   
\end{proof}

\section{The Elliptic Umbilic}

Consider the elliptic umbilic in dimension 2, whose generating
function is
\begin{equation}
\label{ellipticumbilic}
f(y_1,y_2)=y_1^3-y_1y_2^2
\end{equation}
The system (10) in \cite{M}
takes the form
\begin{equation}
\label{dseu}
\left\{ \begin{aligned}
\frac{dy_1}{dt} & =y_1^2-y_2^2-x_1\\
\frac{dy_2}{dt} & =-2y_1y_2-x_2
\end{aligned}
\right.
\end{equation}

\begin{proposition}
\label{nopert}
The bifurcation locus ${\cal B}$ of (\ref{dseu}) consists of three half-lines,
with equations given by $t\rightarrow te^{i\alpha}$, for
$\alpha=0, 2\pi/3, 4\pi/3$, and $t>0$,
and is represented in figure \ref{bifellumb}.1.\\
\end{proposition}
\setlength{\unitlength}{1cm}
\label{bifellumb}
\begin{picture}(4,4)
\thicklines

\qbezier(2,2)(3,2)(4,2)
\qbezier(2,2)(1.5,1.13)(1,0.26)
\qbezier(2,2)(1.5,2.87)(1,3.74)

\thinlines
\put(0,2){\vector(1,0){4}}
\put(2,0){\vector(0,1){4}}
\put(3.7,2.2){$\displaystyle x_1$}
\put(2.2,3.7){$\displaystyle x_2$}

\end{picture}
\begin{center}
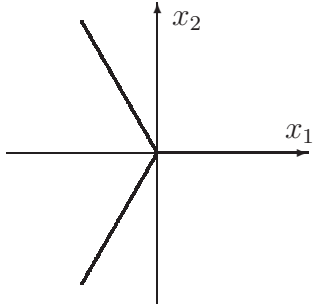
$Fig.~\ref{bifellumb}.1:~The~bifurcation~diagram~of~the~elliptic~umbilic~in
~dimension~2$
\end{center}

\begin{proof}
For $(x_1,x_2)\neq(0,0)$, system (\ref{dseu}) has two critical
points: the saddles $s_1$ and $s_2$
\begin{displaymath}
s_1=\Bigg(\sqrt{\frac{\sqrt{x_1^2+x_2^2}}{2}+x_1},-sgn(x_2)
\sqrt{\frac{\sqrt{x_1^2+x_2^2}}{2}-x_1}\Bigg)
\end{displaymath}
\begin{displaymath}
s_2=\Bigg(-\sqrt{\frac{\sqrt{x_1^2+x_2^2}}{2}+x_1},sgn(x_2)
\sqrt{\frac{\sqrt{x_1^2+x_2^2}}{2}-x_1}\Bigg)
\end{displaymath}
(on the caustic $(x_1,x_2)=(0,0)$, $s_1$ and $s_2$ glue
together in a 2-fold
saddle, also called saddle of multiplicity 2).
Note that $s_1$ and $s_2$ are symmetric with respect to the origin. 
Moreover, let $(\rho,\theta)$ and $(r,\alpha)$ 
be polar
coordinates 
in the $(y_1,y_2)$-plane and $(x_1,x_2)$-plane respectively; suppose $x_2>0$
for simplicity, then the phase 
of the saddles depends on the phase of
the parameter $(x_1,x_2)$ as follows:
$$tg\theta=\frac{y_2}{y_1}=\frac{x_1}{x_2}-\sqrt{\Big(\frac{x_1}{x_2}\Big)^2+1}
=\frac{cos\alpha-1}{sin\alpha}=tg\Big(-\frac{\alpha}{2}\Big)$$
This means that rotating the point $(x_1,x_2)$ clockwise of an angle $\alpha$,
the saddles rotate anti-clockwise of an angle $\alpha/2$.
Observe also that system (\ref{dseu}) can be scaled, in the sense that,
if $y(t)$ is a solution of (\ref{dseu}) corresponding to $x$, then
$\tilde{y}(t)=\frac{1}{\lambda}y(\frac{t}{\lambda})$ is a solution of  
(\ref{dseu})
corresponding to $x/\lambda^2$.
So if a saddle-to-saddle separatrix 
exists for a given $x$, then it exists also for
a positive multiple of $x$. This implies that the bifurcation locus is
formed by rays with source in $(0,0)$.
 
Note also
that if $y(t)$ is a solution of (\ref{dseu}) 
then $-y(t)$ is also a solution
of (\ref{dseu}). In particular, if a saddle-to-saddle separatrix 
exists, being the two
saddles symmetric with respect to the origin, it implies that the 
saddle-to-saddle separatrix 
itself is symmetric with respect to the origin.

In polar coordinates, the generating function is written as
$$f(\rho,\theta)=\rho^3cos\theta(\frac{1}{3}cos^2\theta-sin^2\theta)$$
and (\ref{dseu}) as
\begin{equation}
\label{dseupc}
\left\{ \begin{aligned}
\frac{d\rho}{dt} & =\rho^2cos(3\theta)-rcos(\theta-\alpha)\\
\frac{d\theta}{dt} & =\rho^2sin(3\theta)+rsin(\theta-\alpha)
\end{aligned}
\right.
\end{equation}
We look for solutions through the origin whose image is a straight line, 
imposing the condition $\theta$ 
constant. In order to get a non-constant solution, the
second equation of (\ref{dseupc}) implies that
\begin{equation}
\label{dseupc2}
\left\{ \begin{aligned}
sin(3\theta)=0\\
sin(\theta-\alpha)=0
\end{aligned}
\right.
\end{equation}
These equations are solved by $\theta=k\pi/3$, with $k\in\Z$, and $\alpha=
\theta+h\pi$, with $h\in\Z$. Only for $k$ even we get a
saddle-to-saddle separatrix. 
So ${\cal B}$ contains at least the half-lines
with phase 0, $2\pi/3$ and $4\pi/3$. 
Numerical evidences suggest that no other half-line from the
origin belongs to ${\cal B}$.
\end{proof}

%
%
%

\begin{remark}
\rm
Fukaya in \cite{F2} proposed a conjecture according to which the bifurcation
locus is isotopic to a set of certain integral curves
of the gradient field of the multivalued 
function
\begin{displaymath}
x\mapsto\int_{D^2}u^*\omega
\end{displaymath}
where $u$ is a pseudo-holomorphic disc bounded by the fibre
over $x$ and by the given Lagrangian submanifold. In this particular case the
conjecture is verified.
\end{remark}

The next step is 
to determine the bifurcation diagram when a small perturbation is
added. A first result is lemma 
\ref{far}, which gives information about the structure of the bifurcation
locus outside a disc containing the caustic. 

\begin{lemma}
\label{far}
The bifurcation locus ${\cal B}(\tilde{f})$ of a small perturbation 
$\tilde{f}=f+f'$ of $f$
is, outside a compact subset $D$ containing the support of $f'$,
diffeomorphic to the bifurcation locus ${\cal B}(f)$ of $f$.
\end{lemma}
\begin{proof}
We reason as in proposition (4.10).
We repeat here again the argument.
Let $x_0\in{\cal B}(f)\setminus D$, and
call $\gamma$ the saddle-to-saddle separatrix 
between the saddles $s_1(x_0)$ and 
$s_2(x_0)$ of $\nabla{f_{x_0}}$. 
Consider a 
transversal $\gamma^\perp$ to $\gamma$, identify $\gamma^\perp$ 
with some interval 
$(a,b)$, and set, for every point $x$, 
$h(x)=W^u(s_1(x))\cap\gamma^\perp$ and
$k(x)=W^s(s_2(x))\cap\gamma^\perp$
(note that, since $x\ni D$, $f'_{x}=0$, thus $s_1(x)$ and 
$s_2(x)$ are also saddles of $\nabla\tilde{f}_{x}$). 
This defines
a map $\psi_f:U(x_0)\rightarrow\R$, $\psi_f(x)=h(x)-k(x)$, where $U(x_0)$ is a 
suitably small disc 
around 
$x_0$. The bifurcation locus is the subset $\psi_f^{-1}(0)$. Note that, 
since $U(x_0)\setminus\psi_f^{-1}(0)$ has two connected component 
corresponding to different signs of $\psi_f$, 
the map $\psi_{\tilde{f}}$, defined for a small perturbation
$\tilde{f}$ of $f$, still
attain the value 0 at a certain point
$\tilde{x}$. Moreover, since $\psi$ is a submersion at
$x_0$, it is
transversal to 0, and the point $\tilde{x}$ is unique.
\end{proof}

It remains to study the bifurcation diagram near the caustic. We proceed by 
first 
determining the allowed bifurcation diagrams, 
then theorem (4.30) in \cite{M}
will ensure the existence of perturbations of $f$ exhibiting such diagrams.
As a first step, consider a perturbation given by
a polynomial of degree 2, for example
$f'(y_1,y_2)=\frac{1}{2}(y_1^2+y_2^2)$: 
the generating function is
\begin{equation}
\label{eugfp}
\tilde{f}(y_1,y_2)=\frac{1}{3}y_1^3-2y_1y_2^2+\frac{1}{2}(y_1^2+y_2^2)
\end{equation}
and (\ref{dseu}) becomes
\begin{equation}
\label{dseup}
\left\{ \begin{aligned}
\frac{dy_1}{dt} & =y_1^2-y_2^2+y_1-x_1\\
\frac{dy_2}{dt} & =-2y_1y_2+y_2-x_2
\end{aligned}
\right.
\end{equation}\\
This is the only case where some computation is still feasible.

\begin{proposition}
\label{propos}
The bifurcation locus of (\ref{eugfp}) contains at least three
half-lines departing from each vertex of
the caustic, given by $t\rightarrow te^{i\alpha}$, for
$\alpha=0, 2\pi/3, 4\pi/3$, $t>0$, and represented in figure 
\ref{bifdiaeupol2}.2. 
\end{proposition}
\setlength{\unitlength}{1cm}
\begin{picture}(8,8)
\label{bifdiaeupol2}
\thinlines
\put(0,4){\vector(1,0){8}}
\put(3,0){\vector(0,1){8}}
\put(7.5,4.2){$\displaystyle x_1$}
\put(3.2,7.7){$\displaystyle x_2$}
\qbezier(6,4)(3,4)(2,6)
\qbezier(6,4)(3,4)(2,2)
\qbezier(2,6)(3,4)(2,2)

\thicklines
\qbezier(6,4)(7,4)(8,4)
\qbezier(2,6)(1.5,6.87)(1,7.74)
\qbezier(2,2)(1.5,1.13)(1,0.26)


\end{picture}
\begin{center}
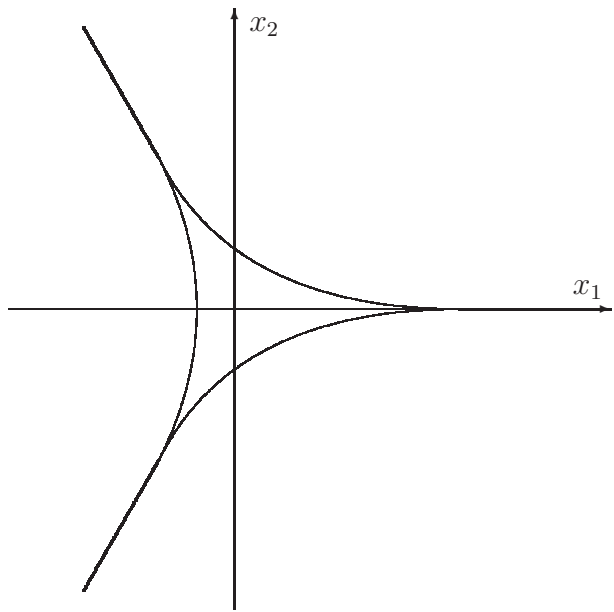
$Fig.~\ref{bifdiaeupol2}.2:~The~bifurcation~diagram~of~a~perturbation~of~the~
elliptic$
$umbilic~by~a~polynomial~of~degree~2$
\end{center}

\begin{proof}
If $x_2=0$ we find the following critical points:
$s_1=\Big(\frac{1}{2},\sqrt{\frac{3}{4}-x_1}\Big)$ and
$s_2=\Big(\frac{1}{2},-\sqrt{\frac{3}{4}-x_1}\Big)$, 
defined for $x_1\leq\frac{3}{4}$, and 
$s_3=\Big(\frac{-1-\sqrt{1+4x_1}}{2},0\Big)$ and 
$n=\Big(\frac{-1+\sqrt{1+4x_1}}{2},0\Big)$, 
defined for $x_1\geq-\frac{1}{4}$.
Linearizing (\ref{dseup}) at these points, we see that
$s_1$, $s_2$ and $s_3$ are saddles, while $n$ is 
an unstable node for $-\frac{1}{4}\leq x_1<\frac{3}{4}$ and a saddle for
$x_1>\frac{3}{4}$. The points $s_1$, $s_2$ and $n$ at $x_1=\frac{3}{4}$,
a vertex of the tricuspoid, 
glue together into a degenerate critical point, which, 
for $x_1>\frac{3}{4}$, turns into a simple (non-degenerate) saddle. Instead,
at $x_1=-\frac{1}{4}$,
$s_3$ and $n$ glue together into a saddle-node, 
which disappears for
$x_1\geq-\frac{1}{4}$.
This implies that, for $(x_1,x_2)$ inside the caustic, system
(\ref{dseup}) has four
critical points, three saddles and an unstable node, while outside the
caustic there are two saddles. On each side of the caustic the node glues
together with one of the three saddles, forming 
a degenerate critical point which disappears outside the caustic.

A gradient line $\gamma_{ns_3}$
from $n$ to $s_3$ can be explicitly computed, 
implying that the half-line
$x_2=0$, $x_1>\frac{3}{4}$, belongs to the
bifurcation locus:
setting $p=\frac{-1+\sqrt{1+4x_1}}{2}$ and $q=\frac{-1-\sqrt{1+4x_1}}{2}$, and
choosing as initial condition a point on the $y_1$-axis between $p$ and $q$,
for example 0, $\gamma_{ns_3}$ is given by
\begin{equation*}
\left\{ \begin{aligned}
y_1(t) & =pq\frac{1-e^{(p-q)t}}{q-pe^{(p-q)t}}\\
y_2(t) & =0
\end{aligned}
\right.
\end{equation*}\\
The existence of the remaining bifurcation half-lines, having the other
vertexes of the caustic as limit point, can be proved in a similar way. 

%
%
%
%
%
\end{proof}

\begin{corollary}
If $f'$ is a generic polynomial of degree 2, up to translations, the
bifurcation diagram of $\tilde{f}=f+f'$ is represented in figure
\ref{bifdiaeupol2}.3.
\end{corollary}
\begin{proof}
As done in subsection (3.1) of \cite{M},
we can reduce to the hypothesis of proposition
\ref{propos} as a consequence of a suitable translation.
\end{proof}

Before analysing a generic small perturbation of the elliptic umbilic
(\ref{ellipticumbilic}), consider 
a generic small perturbations of (\ref{eugfp}).

\begin{proposition}
\label{propos2}
The bifurcation locus of a small perturbation of (\ref{eugfp}),
shown in figure \ref{bifpertpol2}.3, has the following features:\\
1. outside a compact subset containing the caustic $K$, there exist 
three bifurcation lines (as for
the unperturbed elliptic umbilic);\\
2. generically, these half-lines intersect $K$ at a point of 
one of its sides near a vertex and have as extreme a fold point on
the opposite side of $K$, near the same vertex;\\
3. as already described for the cusp in proposition \ref{cuspbifu}, 
inside $K$ also
immersed $S^1$'s or segments, with extremes on the same side of $K$, 
may appear as components of the bifurcation locus;\\ 
\end{proposition}

\setlength{\unitlength}{1cm}
\begin{picture}(8,8)
\label{bifpertpol2}
\thinlines
\put(0,4){\vector(1,0){8}}
\put(3,0){\vector(0,1){8}}
\put(7.5,4.2){$\displaystyle x_1$}
\put(3.2,7.7){$\displaystyle x_2$}
\qbezier(6,4)(3,4)(2,6)
\qbezier(6,4)(3,4)(2,2)
\qbezier(2,6)(3,4)(2,2)

\thicklines
\qbezier(8,3.8)(7.9,3.8)(7.4,3.7)
\qbezier(7.4,3.7)(7.1,3.5)(6,3.5)
\qbezier(4.7,4.1)(5,3.5)(6,3.5)

\qbezier(1.7,5.5)(1.25,7.6)(1,7.74)
\qbezier(1.8,5.1)(1.7,5.2)(1.7,5.5)
\qbezier(2.5,5.2)(2.1,5)(1.8,5.1)

\qbezier(2.5,2.2)(1.25,0.7)(1,0.26)
\qbezier(3,3)(2.7,2.3)(2.5,2.2)
\qbezier(3,3)(3,3.2)(2.8,3.2)
\qbezier(2.3,2.7)(2.5,3.2)(2.8,3.2)

\qbezier(3.15,4.65)(3,4)(2.7,5)
\qbezier(3.5,4.1)(3.1,4.2)(3.5,4.4)
\qbezier(3.5,4.1)(3.9,4.1)(3.5,4.4)

\end{picture}
\begin{center}
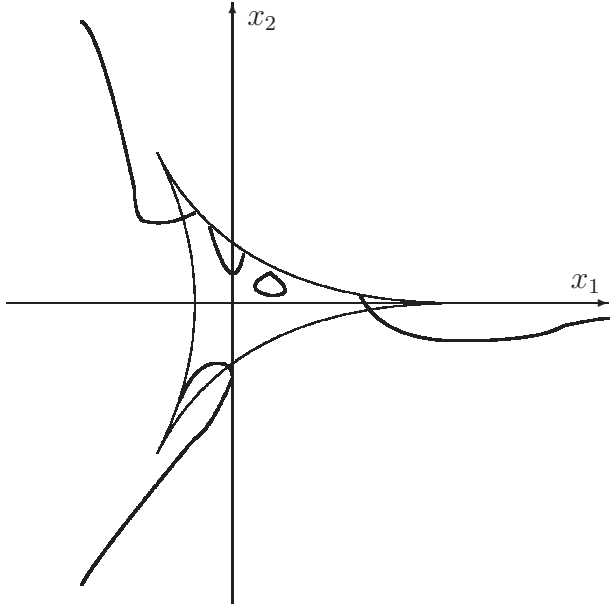
$Fig.~\ref{bifpertpol2}.3:~The~bifurcation~diagram~of~a~small~
perturbation~of~a$
$perturbation~by~a~polynomial~of~degree~2~of~the~
elliptic~umbilic$
\end{center}

\begin{proof}
The structure of the bifurcation locus far from $K$, as described
at point 1, is a 
direct consequence of
lemma \ref{far}. This, together with proposition (4.10) in \cite{M},
implies that the 
bifurcation locus contains half-lines with
extreme on the caustic.
If $x$ is a cusp (a vertex) of $K$, then $\nabla\tilde{f}_x$ 
has a saddle $s_i$ and
a degenerate critical point $ns_js_k$. Suppose a bifurcation line 
${\cal B}$ has
$x$ as limit point, then, for all $t\in{\cal B}\cap N$, where $N$ is
a neighbourhood of the vertex $x$, the field $\nabla\tilde{f}_t$ exhibits
a saddle-to-saddle
separatrix $\gamma_t$, and, at the limit point $x$, a gradient line 
$\gamma_v$ exists
between $s_i$ and $ns_js_k$. Since $\gamma_t$ is unstable, also $\gamma_v$
is unstable. Thus, generically, the limit point of a bifurcation line is
a fold of $K$.
That a bifurcation half-line {\cal B} intersects $K$ on a side 
$l_i$ and has a limit point on one of the opposite sides $l_j$,
is a consequence of the fact that the exceptional gradient line 
$\gamma_{s_js_k}$, exhibited by $\nabla\tilde{f}_x$ for $x\in{\cal B}$, 
breaks when one of the points $s_j$ or $s_k$ glues
together with the node $n$ in a saddle-node, but this happens just when
$x$ belongs to one of the opposite sides to $l_i$.
The argument to prove point 3 is the same given for the cusp in proposition
\ref{cuspbifu}. 
The proposition now follows from theorem (4.30) in \cite{M}.
\end{proof}  

The final step is to determine the bifurcation locus of a generic 
perturbation of (\ref{ellipticumbilic}). 
The behaviour of the bifurcation locus 
outside the caustic is determined by lemma
\ref{far}. The idea to study how the bifurcation locus looks inside
the caustic $K$ is as follows: we consider diagrams representing all the
possible mutual positions of three bifurcation half-lines inside
$K$, having an extreme on a certain side of $K$ and 
intersecting further $K$ at
a point of the remaining sides, first assuming that such
lines do not intersect, and we study which among such diagrams are allowed;
then we do the same assuming bifurcation lines can intersect;
finally, theorem (4.30) in \cite{M}
ensures
the existence of a function $\tilde{f}$
giving rise to such bifurcation diagrams.

As before we denote by $n$ and $s_i$, $i=1,2,3$, the node and the
saddles of $\nabla\tilde{f}_x$, for $x$ lying inside the caustic.   
Let $l_i$, for $i=1,2,3$, be the side of the caustic where the
saddle $s_i$ glues together with $n$. 
Observe that if a bifurcation half-line ${\cal B}$ intersects the side $l_i$
of the caustic and has its extreme on the opposite
side $l_j$, it means that, for $x\in{\cal B}$, 
$\nabla\tilde{f}_x$ exhibits in its phase portrait the
saddle-to-saddle separatrix $\gamma_{s_js_k}$, 
from $s_j$ to $s_k$, with $k\in\{1,2,3\}\setminus\{i,j\}$: in other
words, ${\cal B}\subset{\cal B}_{jk}$. Assume 
that bifurcation lines do not intersect. 
In this case, 
all the possible diagrams are drawn in figure \ref{allowedandnot}.4.\\
\setlength{\unitlength}{1cm}
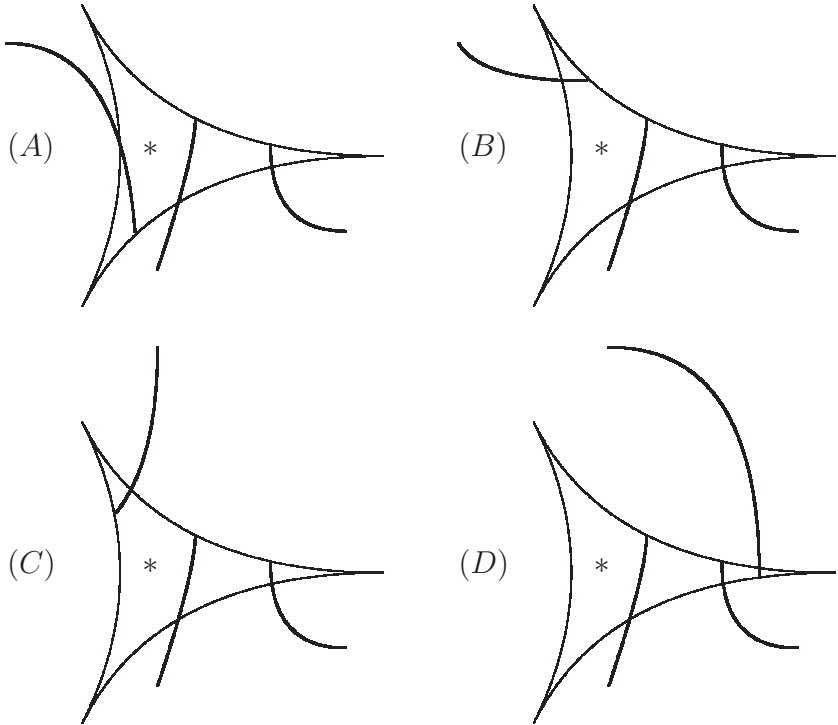
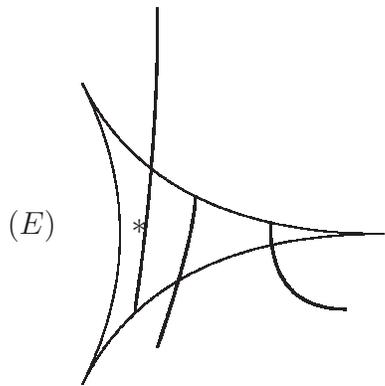
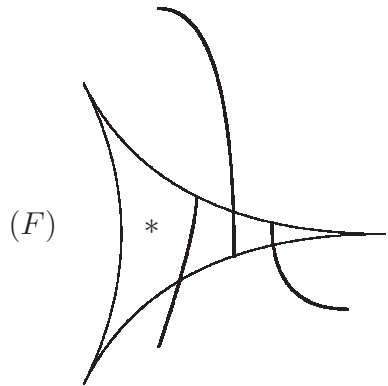
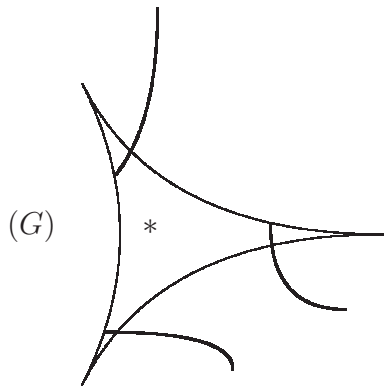
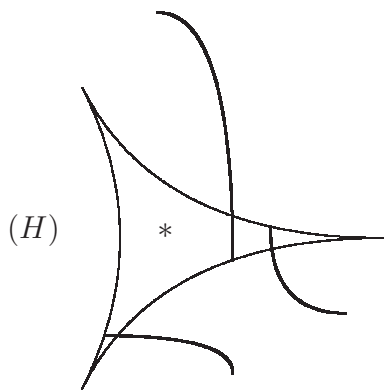
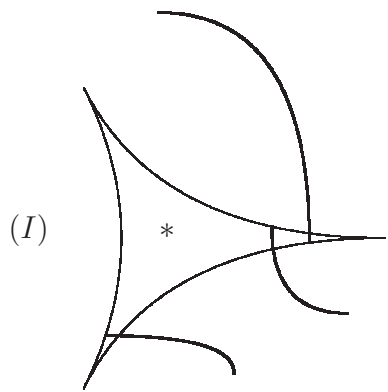
\label{allowedandnot}
\begin{picture}(12,5)

\thinlines
\qbezier(5,2)(2,2)(1,4)
\qbezier(5,2)(2,2)(1,0)
\qbezier(1,4)(2,2)(1,0)

\thicklines
\qbezier(4.5,1)(3.5,1)(3.5,2.15)
\qbezier(2,0.5)(2.5,2)(2.5,2.5)
\qbezier(0,3.5)(1.5,3.5)(1.7,1)

\thinlines
\qbezier(11,2)(8,2)(7,4)
\qbezier(11,2)(8,2)(7,0)
\qbezier(7,4)(8,2)(7,0)

\thicklines
\qbezier(10.5,1)(9.5,1)(9.5,2.15)
\qbezier(8,0.5)(8.5,2)(8.5,2.5)
\qbezier(6,3.5)(6.3,3)(7.7,3)

\put(0,2){$(A)$}
\put(6,2){$(B)$}

\put(1.8,2){$\displaystyle *$}
\put(7.8,2){$\displaystyle *$}

\end{picture}\\\\
\setlength{\unitlength}{1cm}
\begin{picture}(12,5)

\thinlines
\qbezier(5,2)(2,2)(1,4)
\qbezier(5,2)(2,2)(1,0)
\qbezier(1,4)(2,2)(1,0)

\thicklines
\qbezier(4.5,1)(3.5,1)(3.5,2.15)
\qbezier(2,0.5)(2.5,2)(2.5,2.5)
\qbezier(2,5)(2,3.5)(1.45,2.8)

\thinlines
\qbezier(11,2)(8,2)(7,4)
\qbezier(11,2)(8,2)(7,0)
\qbezier(7,4)(8,2)(7,0)

\thicklines
\qbezier(10.5,1)(9.5,1)(9.5,2.15)
\qbezier(8,0.5)(8.5,2)(8.5,2.5)
\qbezier(8,5)(10,5)(10,1.95)

\put(0,2){$(C)$}
\put(6,2){$(D)$}

\put(1.8,2){$\displaystyle *$}
\put(7.8,2){$\displaystyle *$}

\end{picture}\\\\
\setlength{\unitlength}{1cm}
\begin{picture}(12,5)

\thinlines
\qbezier(5,2)(2,2)(1,4)
\qbezier(5,2)(2,2)(1,0)
\qbezier(1,4)(2,2)(1,0)

\thicklines
\qbezier(4.5,1)(3.5,1)(3.5,2.15)
\qbezier(2,0.5)(2.5,2)(2.5,2.5)
\qbezier(2,5)(2,3)(1.7,1)

\thinlines
\qbezier(11,2)(8,2)(7,4)
\qbezier(11,2)(8,2)(7,0)
\qbezier(7,4)(8,2)(7,0)

\thicklines
\qbezier(10.5,1)(9.5,1)(9.5,2.15)
\qbezier(8,0.5)(8.5,2)(8.5,2.5)
\qbezier(8,5)(9,5)(9,1.7)

\put(0,2){$(E)$}
\put(6,2){$(F)$}

\put(1.65,2){$\displaystyle *$}
\put(7.8,2){$\displaystyle *$}

\end{picture}\\\\
\setlength{\unitlength}{1cm}
\begin{picture}(12,5.5)

\thinlines
\qbezier(5,2)(2,2)(1,4)
\qbezier(5,2)(2,2)(1,0)
\qbezier(1,4)(2,2)(1,0)

\thicklines
\qbezier(4.5,1)(3.5,1)(3.5,2.15)
\qbezier(3,0.2)(3,.7)(1.3,0.7)
\qbezier(2,5)(2,3.5)(1.45,2.8)

\put(0,2){$(G)$}

\put(1.8,2){$\displaystyle *$}

\end{picture}\\\\
\setlength{\unitlength}{1cm}
\begin{picture}(12,5)

\thinlines
\qbezier(5,2)(2,2)(1,4)
\qbezier(5,2)(2,2)(1,0)
\qbezier(1,4)(2,2)(1,0)

\thicklines
\qbezier(4.5,1)(3.5,1)(3.5,2.15)
\qbezier(3,0.2)(3,.7)(1.3,0.7)
\qbezier(2,5)(3,5)(3,1.7)

\thinlines
\qbezier(11,2)(8,2)(7,4)
\qbezier(11,2)(8,2)(7,0)
\qbezier(7,4)(8,2)(7,0)

\thicklines
\qbezier(10.5,1)(9.5,1)(9.5,2.15)
\qbezier(9,0.2)(9,.7)(7.3,0.7)
\qbezier(8,5)(10,5)(10,1.95)

\put(0,2){$(H)$}
\put(6,2){$(I)$}

\put(2,2){$\displaystyle *$}
\put(8,2){$\displaystyle *$}

\end{picture}\\\\
\setlength{\unitlength}{1cm}
\begin{picture}(12,5)

\thinlines
\qbezier(5,2)(2,2)(1,4)
\qbezier(5,2)(2,2)(1,0)
\qbezier(1,4)(2,2)(1,0)

\thicklines
\qbezier(4.5,1)(3.5,1)(3.5,2.15)
\qbezier(0,3.5)(0.3,3)(1.7,3)
\qbezier(2,5)(3,5)(3,1.7)

\thinlines
\qbezier(11,2)(8,2)(7,4)
\qbezier(11,2)(8,2)(7,0)
\qbezier(7,4)(8,2)(7,0)

\thicklines
\qbezier(10.5,1)(9.5,1)(9.5,2.15)
\qbezier(6,3.5)(6.3,3)(7.7,3)
\qbezier(8,5)(10,5)(10,1.95)

\put(2,2){$\displaystyle *$}
\put(8,2){$\displaystyle *$}

\put(0,2){$(J)$}
\put(6,2){$(K)$}

\end{picture}\\\\
\setlength{\unitlength}{1cm}
\begin{picture}(12,5)

\thinlines
\qbezier(5,2)(2,2)(1,4)
\qbezier(5,2)(2,2)(1,0)
\qbezier(1,4)(2,2)(1,0)

\thicklines
\qbezier(4.5,1)(3.5,1)(3.5,2.15)
\qbezier(0,1.5)(1.5,1.5)(1.7,1)
\qbezier(2,5)(2,3.5)(1.45,2.8)

\put(2,2){$\displaystyle *$}

\put(0,2){$(L)$}

\end{picture}\\\\
\setlength{\unitlength}{1cm}
\begin{picture}(12,5)

\thinlines
\qbezier(5,2)(2,2)(1,4)
\qbezier(5,2)(2,2)(1,0)
\qbezier(1,4)(2,2)(1,0)

\thicklines
\qbezier(4.5,1)(3.5,1)(3.5,2.15)
\qbezier(3,0.5)(2.5,2)(2.5,2.5)
\qbezier(2,0.2)(2,.7)(1.3,0.7)

\thinlines
\qbezier(11,2)(8,2)(7,4)
\qbezier(11,2)(8,2)(7,0)
\qbezier(7,4)(8,2)(7,0)

\thicklines
\qbezier(10.5,1)(9.5,1)(9.5,2.15)
\qbezier(9,0.5)(8.5,2)(8.5,2.5)
\qbezier(8,0.2)(8,.7)(7.7,3)

\put(0,2){$(M)$}
\put(6,2){$(N)$}

\put(1.8,2){$\displaystyle *$}
\put(7.55,2){$\displaystyle *$}

%
%
%
%
%
%
%
\end{picture}\\
\begin{center}
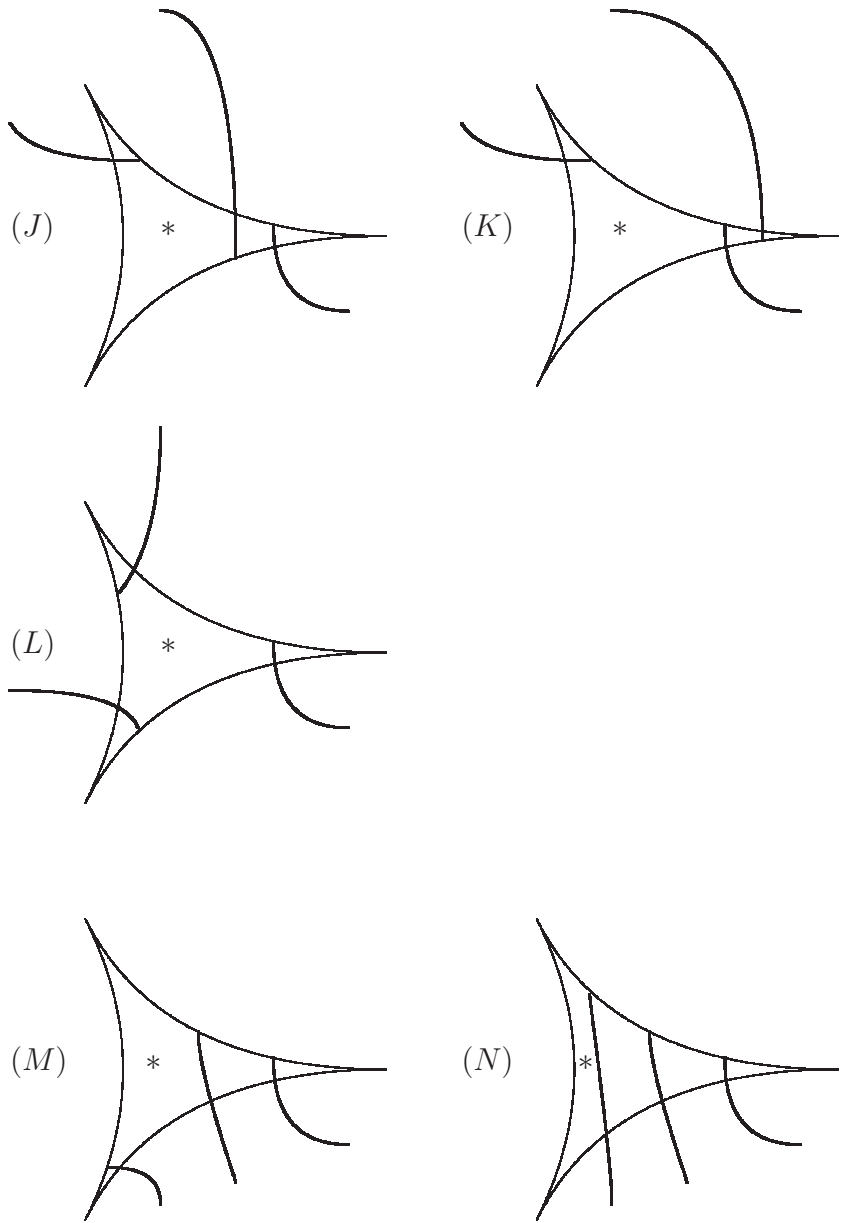
$Fig.~\ref{allowedandnot}.4:~Possible~mutual~positions~of~non-intersecting$
$bifurcation~lines
~inside~the~caustic$
\end{center}

\begin{lemma}
\label{lemdiagr}
Among diagrams in figure \ref{allowedandnot}.4, only (A), (B), (C),
(D), (G), (L), (M), (N) are allowed.
\end{lemma}

\begin{proof}
In all diagrams, the subset $\ast$ is the one bounded
by all the sides of the caustic. For
every $x\in\ast$, the phase portrait of $\nabla\tilde{f}_x$
contains all the gradient lines $\gamma_{ns_i}$, $i=1,2,3$: indeed, as
already explained in the proof of proposition \ref{cuspbifu}, if 
$\gamma_{ns_i}$ were missing, then, for $x$ belonging to 
the side $l_i$ of the caustic
bounding $\ast$, $\nabla\tilde{f}_x$ would exhibit
two gradient lines between the $ns_i$ and
one of the remaining saddles, implying a contraddiction.
The lemma follows now from 
proposition (4.19) in \cite{M}
with considerations as those in subsection (4.5) of \cite{M}.
\end{proof}

We study now the possibility of intersection between bifurcation lines.
The possible intersections, depending on the positions of bifurcation
lines, are listed in figure \ref{cases}.5.\\
\setlength{\unitlength}{1cm}
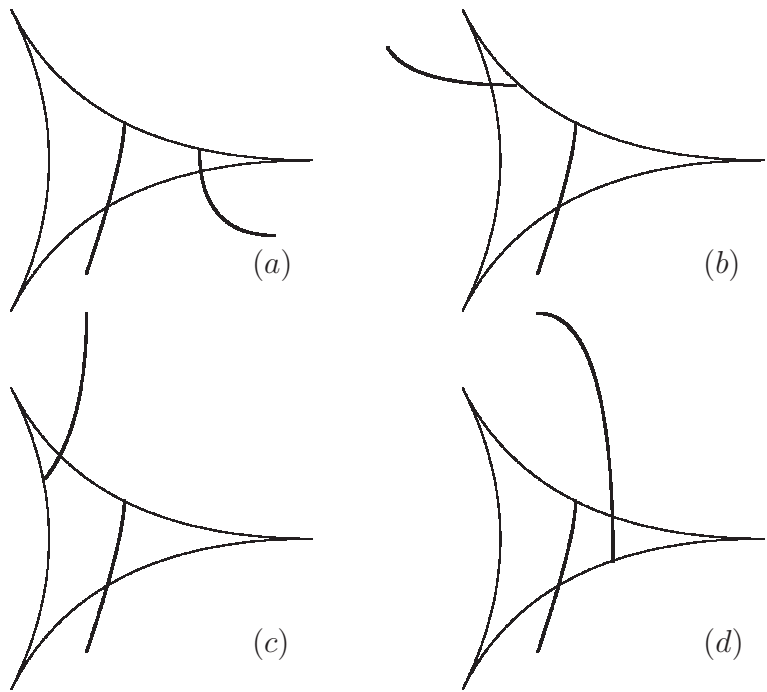
\begin{picture}(12,5)
\label{cases}

\thinlines
\qbezier(5,2)(2,2)(1,4)
\qbezier(5,2)(2,2)(1,0)
\qbezier(1,4)(2,2)(1,0)

\thicklines
\qbezier(4.5,1)(3.5,1)(3.5,2.15)
\qbezier(2,0.5)(2.5,2)(2.5,2.5)

\thinlines
\qbezier(11,2)(8,2)(7,4)
\qbezier(11,2)(8,2)(7,0)
\qbezier(7,4)(8,2)(7,0)

\thicklines
\qbezier(8,0.5)(8.5,2)(8.5,2.5)
\qbezier(6,3.5)(6.3,3)(7.7,3)

\put(4.2,.5){$(a)$}
\put(10.2,.5){$(b)$}

\end{picture}\\
\setlength{\unitlength}{1cm}
\begin{picture}(12,5)

\thinlines
\qbezier(5,2)(2,2)(1,4)
\qbezier(5,2)(2,2)(1,0)
\qbezier(1,4)(2,2)(1,0)

\thicklines
\qbezier(2,0.5)(2.5,2)(2.5,2.5)
\qbezier(2,5)(2,3.5)(1.45,2.8)

\thinlines
\qbezier(11,2)(8,2)(7,4)
\qbezier(11,2)(8,2)(7,0)
\qbezier(7,4)(8,2)(7,0)

\thicklines
\qbezier(8,0.5)(8.5,2)(8.5,2.5)
\qbezier(8,5)(9,5)(9,1.7)

\put(4.2,.5){$(c)$}
\put(10.2,.5){$(d)$}

\end{picture}\\

\setlength{\unitlength}{1cm}
\begin{picture}(12,5)

\thinlines
\qbezier(5,2)(2,2)(1,4)
\qbezier(5,2)(2,2)(1,0)
\qbezier(1,4)(2,2)(1,0)

\thicklines
\qbezier(4.5,1)(3.5,1)(3.5,2.15)
\qbezier(2,0.2)(2,.7)(1.3,0.7)

\put(4.2,.5){$(e)$}


\end{picture}
\begin{center}
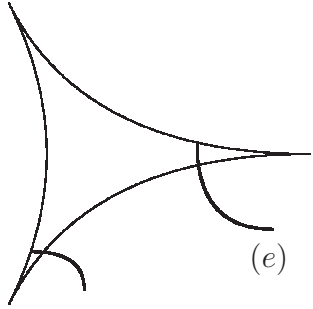
$Fig.~\ref{cases}.5:~Possible~cases~of~intersection~of~bifurcation~lines$
\end{center}
In all these diagrams, the third bifurcation half-line, which is not shown,
has a position
such that the resulting bifurcation diagram is among those allowed by lemma
\ref{lemdiagr}. In what follows, $l_1$ denotes the left side of the
caustic, $l_2$ its upper side and $l_3$ its lower side.

Consider (a) (see figure \ref{casoa}.6).\\
\setlength{\unitlength}{1cm}
\label{casoa}
\begin{picture}(12,4.5)

\thinlines
\qbezier(5,2)(2,2)(1,4)
\qbezier(5,2)(2,2)(1,0)
\qbezier(1,4)(2,2)(1,0)

\thicklines
\qbezier(4.5,1)(3.5,1)(3.5,2.15)
\qbezier(2,0.5)(2.5,2)(2.5,2.5)

\thinlines
\qbezier(11,2)(8,2)(7,4)
\qbezier(11,2)(8,2)(7,0)
\qbezier(7,4)(8,2)(7,0)

\thicklines
\qbezier(10.5,1)(9.5,1)(8.5,2.5)
\qbezier(8,0.5)(8.5,2)(9,2.25)

\put(1.8,1.9){$\alpha$}
\put(2.9,1.9){$\beta$}
\put(3.7,1.9){$\gamma$}

\put(7.8,1.9){$\alpha$}
\put(8.7,1.6){$\beta$}
\put(9.4,1.9){$\gamma$}
\put(8.7,2.2){$\delta$}


\end{picture}
\begin{center}
$Fig.~\ref{casoa}.6:~Intersection~of~bifurcation~lines:~case~(a)$
\end{center}

\begin{proposition}
\label{propcasoa}
The intersection of bifurcation lines of $(a)$ gives an allowed
bifurcation diagram provided
$\nabla\tilde{f}_x$ exhibits, at points of both bifurcation lines,
a saddle-to-saddle separatrix obtained by
joining the
same pair of separatrices.
\end{proposition}
\begin{proof}
The result follows from lemma (4.18) in \cite{M}.
At points of
bifurcation
lines the gradient line $\gamma_{s_2s_1}$ appears in the phase portrait
of $\nabla\tilde{f}_x$,
however this can occur in three ways:\\
1) the non-generic gradient line $\gamma_{s_2s_1}$ in both bifurcations
is obtained by joining the same pair of separatrices; 
the phase portrait of $\nabla\tilde{f}_x$ for $x\in\alpha$ is the same
as the one for $x\in\gamma$ and it contains the gradient line $\gamma_{ns_1}$,
which disappears for $x\in\beta$ (see figure \ref{primabif}.7).\\
\setlength{\unitlength}{1cm}
\begin{picture}(14,4.5)
\label{primabif}

\thinlines
\put(1,2){\circle*{.1}}
\put(2.7,2.7){\circle*{.1}}
\put(2.7,1.3){\circle*{.1}}
\put(2,2){\circle*{.1}}

\scriptsize
\put(1.1,1.7){$\displaystyle s_1$}
\put(2.2,1.9){$\displaystyle n$}
\put(2.9,2.7){$\displaystyle s_2$}
\put(2.9,1.2){$\displaystyle s_3$}
\normalsize
\put(0,.5){$\displaystyle (\alpha)$}

\qbezier(2,2)(1.5,2)(1,2)
\qbezier(2,2)(2.5,2.5)(2.7,2.7)
\qbezier(2,2)(2.5,1.5)(2.7,1.3)

\qbezier(1,2)(1,3)(1,3.5)
\qbezier(1,2)(1,1)(1,.5)
\qbezier(1,2)(.5,2)(0,2)

\qbezier(2.7,2.7)(3,3)(3.5,3.5)
\qbezier(2.7,2.7)(2.4,3)(2,3.4)
\qbezier(2.7,2.7)(3,2.4)(3.5,2.2)

\qbezier(2.7,1.3)(3,1)(3.5,.6)
\qbezier(2.7,1.3)(3,1.6)(3.5,1.8)
\qbezier(2.7,1.3)(2.4,1)(2,.6)

\thinlines
\put(7,2){\circle*{.1}}
\put(8.7,2.7){\circle*{.1}}
\put(8.7,1.3){\circle*{.1}}
\put(8,2){\circle*{.1}}

\scriptsize
\put(7.1,1.8){$\displaystyle s_1$}
\put(8.2,1.9){$\displaystyle n$}
\put(8.9,2.7){$\displaystyle s_2$}
\put(8.9,1.2){$\displaystyle s_3$}
\normalsize
\put(6,.5){$\displaystyle (\beta)$}

\qbezier(8,3.2)(7.5,2.5)(7,2)
\qbezier(8,2)(8.5,2.5)(8.7,2.7)
\qbezier(8,2)(8.5,1.5)(8.7,1.3)

\qbezier(7,2)(7,3)(7,3.5)
\qbezier(7,2)(7,1)(7,.5)
\qbezier(7,2)(6.5,2)(6,2)

\qbezier(8.7,2.7)(9,3)(9.5,3.5)
\qbezier(8.7,2.7)(7.3,2.2)(7.3,.7)
\qbezier(8.7,2.7)(9,2.4)(9.5,2.2)

\qbezier(8.7,1.3)(9,1)(9.5,.6)
\qbezier(8.7,1.3)(9,1.6)(9.5,1.8)
\qbezier(8.7,1.3)(8.4,1)(8,.6)

\put(4.5,2){\vector(1,0){1}}
\put(10.5,2){\vector(1,0){1}}
\put(4.7,2.3){$1^{st}$}
\put(10.7,2.3){$2^{nd}$}

\end{picture}\\
\setlength{\unitlength}{1cm}
\begin{picture}(14,5)

\thinlines
\put(1,2){\circle*{.1}}
\put(2.7,2.7){\circle*{.1}}
\put(2.7,1.3){\circle*{.1}}
\put(2,2){\circle*{.1}}

\scriptsize
\put(1.1,1.7){$\displaystyle s_1$}
\put(2.2,1.9){$\displaystyle n$}
\put(2.9,2.7){$\displaystyle s_2$}
\put(2.9,1.2){$\displaystyle s_3$}
\normalsize
\put(0,.5){$\displaystyle (\gamma)$}

\qbezier(2,2)(1.5,2)(1,2)
\qbezier(2,2)(2.5,2.5)(2.7,2.7)
\qbezier(2,2)(2.5,1.5)(2.7,1.3)

\qbezier(1,2)(1,3)(1,3.5)
\qbezier(1,2)(1,1)(1,.5)
\qbezier(1,2)(.5,2)(0,2)

\qbezier(2.7,2.7)(3,3)(3.5,3.5)
\qbezier(2.7,2.7)(2.4,3)(2,3.4)
\qbezier(2.7,2.7)(3,2.4)(3.5,2.2)

\qbezier(2.7,1.3)(3,1)(3.5,.6)
\qbezier(2.7,1.3)(3,1.6)(3.5,1.8)
\qbezier(2.7,1.3)(2.4,1)(2,.6)

\end{picture}
\begin{center}
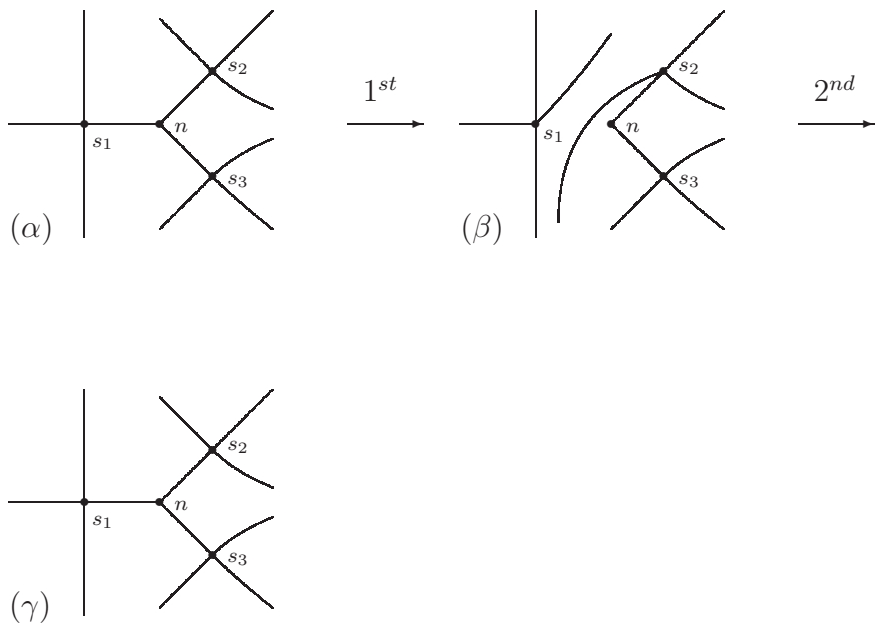
$Fig.~\ref{primabif}.7:~The~bifurcation~from~\alpha~to~\beta~and~from~\beta~
to~\gamma~respectively:$
$1^{st}~possibility$
\end{center}
2) the non-generic gradient line $\gamma_{s_2s_1}$
is obtained by joining different pairs of separatrices; 
the phase portrait of $\nabla\tilde{f}_x$ for $x\in\alpha$ differs
from the one for $x\in\gamma$ (see figure \ref{secondabif}.8).\\
\setlength{\unitlength}{1cm}
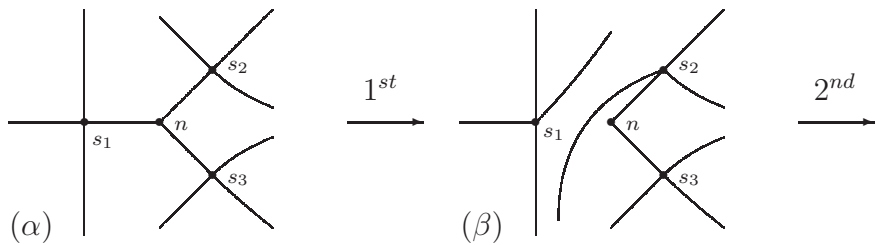
\begin{picture}(14,4.5)
\label{secondabif}

\thinlines
\put(1,2){\circle*{.1}}
\put(2.7,2.7){\circle*{.1}}
\put(2.7,1.3){\circle*{.1}}
\put(2,2){\circle*{.1}}

\scriptsize
\put(1.1,1.7){$\displaystyle s_1$}
\put(2.2,1.9){$\displaystyle n$}
\put(2.9,2.7){$\displaystyle s_2$}
\put(2.9,1.2){$\displaystyle s_3$}
\normalsize
\put(0,.5){$\displaystyle (\alpha)$}

\qbezier(2,2)(1.5,2)(1,2)
\qbezier(2,2)(2.5,2.5)(2.7,2.7)
\qbezier(2,2)(2.5,1.5)(2.7,1.3)

\qbezier(1,2)(1,3)(1,3.5)
\qbezier(1,2)(1,1)(1,.5)
\qbezier(1,2)(.5,2)(0,2)

\qbezier(2.7,2.7)(3,3)(3.5,3.5)
\qbezier(2.7,2.7)(2.4,3)(2,3.4)
\qbezier(2.7,2.7)(3,2.4)(3.5,2.2)

\qbezier(2.7,1.3)(3,1)(3.5,.6)
\qbezier(2.7,1.3)(3,1.6)(3.5,1.8)
\qbezier(2.7,1.3)(2.4,1)(2,.6)

\thinlines
\put(7,2){\circle*{.1}}
\put(8.7,2.7){\circle*{.1}}
\put(8.7,1.3){\circle*{.1}}
\put(8,2){\circle*{.1}}

\scriptsize
\put(7.1,1.8){$\displaystyle s_1$}
\put(8.2,1.9){$\displaystyle n$}
\put(8.9,2.7){$\displaystyle s_2$}
\put(8.9,1.2){$\displaystyle s_3$}
\normalsize
\put(6,.5){$\displaystyle (\beta)$}

\qbezier(8,3.2)(7.5,2.5)(7,2)
\qbezier(8,2)(8.5,2.5)(8.7,2.7)
\qbezier(8,2)(8.5,1.5)(8.7,1.3)

\qbezier(7,2)(7,3)(7,3.5)
\qbezier(7,2)(7,1)(7,.5)
\qbezier(7,2)(6.5,2)(6,2)

\qbezier(8.7,2.7)(9,3)(9.5,3.5)
\qbezier(8.7,2.7)(7.3,2.2)(7.3,.7)
\qbezier(8.7,2.7)(9,2.4)(9.5,2.2)

\qbezier(8.7,1.3)(9,1)(9.5,.6)
\qbezier(8.7,1.3)(9,1.6)(9.5,1.8)
\qbezier(8.7,1.3)(8.4,1)(8,.6)

\put(4.5,2){\vector(1,0){1}}
\put(10.5,2){\vector(1,0){1}}
\put(4.7,2.3){$1^{st}$}
\put(10.7,2.3){$2^{nd}$}

\end{picture}\\
\setlength{\unitlength}{1cm}
\begin{picture}(14,5)

\thinlines
\put(1,2){\circle*{.1}}
\put(2.7,2.7){\circle*{.1}}
\put(2.7,1.3){\circle*{.1}}
\put(2,2){\circle*{.1}}

\scriptsize
\put(.9,1.7){$\displaystyle s_1$}
\put(2.2,1.9){$\displaystyle n$}
\put(2.9,2.7){$\displaystyle s_2$}
\put(2.9,1.2){$\displaystyle s_3$}
\normalsize
\put(0,.5){$\displaystyle (\gamma)$}

\qbezier(1.2,3.5)(1.2,2)(1,2)
\qbezier(2,2)(2.5,2.5)(2.7,2.7)
\qbezier(2,2)(2.5,1.5)(2.7,1.3)

\qbezier(1,2)(1,3)(1,3.5)
\qbezier(1,2)(1.4,2)(1.4,3.4)
\qbezier(1,2)(.5,2)(0,2)

\qbezier(2.7,2.7)(1,2.3)(1,.8)
\qbezier(2.7,2.7)(1.5,2,2)(1.5,.7)
\qbezier(2.7,2.7)(3,2.4)(3.5,2.2)

\qbezier(2.7,1.3)(3,1)(3.5,.6)
\qbezier(2.7,1.3)(3,1.6)(3.5,1.8)
\qbezier(2.7,1.3)(2.4,1)(2,.6)



\end{picture}
\begin{center}
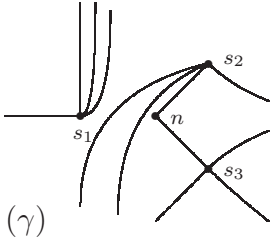
$Fig.~\ref{primabif}.8:~The~bifurcation~from~\alpha~to~\beta~and~from~\beta~
to~\gamma~respectively:$
$2^{nd}~possibility$
\end{center}
3) the non-generic gradient line $\gamma_{s_2s_1}$ in both bifurcations
is obtained by joining the same pair of separatrices; 
the phase portrait of $\nabla\tilde{f}_x$ for $x\in\alpha$ is the same
as the one 
for $x\in\gamma$ and it never contains the gradient line $\gamma_{ns_1}$
(see figure \ref{terzabif}.9).\\
\setlength{\unitlength}{1cm}
\begin{picture}(14,4.5)
\label{terzabif}

\thinlines
\put(1,2){\circle*{.1}}
\put(2.7,2.7){\circle*{.1}}
\put(2.7,1.3){\circle*{.1}}
\put(2,2){\circle*{.1}}

\scriptsize
\put(1.1,1.8){$\displaystyle s_1$}
\put(2.2,1.9){$\displaystyle n$}
\put(2.9,2.7){$\displaystyle s_2$}
\put(2.9,1.2){$\displaystyle s_3$}
\normalsize
\put(0,.5){$\displaystyle (\alpha)$}

\qbezier(2,2)(1.5,2)(1,2)
\qbezier(2,2)(2.5,2.5)(2.7,2.7)
\qbezier(2,2)(2.5,1.5)(2.7,1.3)

\qbezier(1,2)(1,3)(1,3.5)
\qbezier(1,2)(1,1)(1,.5)
\qbezier(1,2)(.5,2)(0,2)

\qbezier(2.7,2.7)(3,3)(3.5,3.5)
\qbezier(2.7,2.7)(2.4,3)(2,3.4)
\qbezier(2.7,2.7)(3,2.4)(3.5,2.2)

\qbezier(2.7,1.3)(3,1)(3.5,.6)
\qbezier(2.7,1.3)(3,1.6)(3.5,1.8)
\qbezier(2.7,1.3)(2.4,1)(2,.6)

\thinlines
\put(7,2){\circle*{.1}}
\put(8.7,2.7){\circle*{.1}}
\put(8.7,1.3){\circle*{.1}}
\put(8,2){\circle*{.1}}

\scriptsize
\put(6.9,1.7){$\displaystyle s_1$}
\put(8.2,1.9){$\displaystyle n$}
\put(8.9,2.7){$\displaystyle s_2$}
\put(8.9,1.2){$\displaystyle s_3$}
\normalsize
\put(6,.5){$\displaystyle (\beta)$}

\qbezier(7.2,3.5)(7.2,2)(7,2)
\qbezier(8,2)(8.5,2.5)(8.7,2.7)
\qbezier(8,2)(8.5,1.5)(8.7,1.3)

\qbezier(7,2)(7,3)(7,3.5)
\qbezier(7,2)(7.4,2)(7.4,3.4)
\qbezier(7,2)(6.5,2)(6,2)

\qbezier(8.7,2.7)(7,2.3)(7,.8)
\qbezier(8.7,2.7)(7.5,2,2)(7.5,.7)
\qbezier(8.7,2.7)(9,2.4)(9.5,2.2)

\qbezier(8.7,1.3)(9,1)(9.5,.6)
\qbezier(8.7,1.3)(9,1.6)(9.5,1.8)
\qbezier(8.7,1.3)(8.4,1)(8,.6)

\end{picture}\\
\setlength{\unitlength}{1cm}
\begin{picture}(14,5)

\thinlines
\put(1,2){\circle*{.1}}
\put(2.7,2.7){\circle*{.1}}
\put(2.7,1.3){\circle*{.1}}
\put(2,2){\circle*{.1}}

\scriptsize
\put(1.1,1.8){$\displaystyle s_1$}
\put(2.2,1.9){$\displaystyle n$}
\put(2.9,2.6){$\displaystyle s_2$}
\put(2.9,1.2){$\displaystyle s_3$}
\normalsize
\put(0,.5){$\displaystyle (\gamma)$}

\qbezier(2,3.2)(1.5,2.5)(1,2)
\qbezier(2,2)(2.5,2.5)(2.7,2.7)
\qbezier(2,2)(2.5,1.5)(2.7,1.3)

\qbezier(1,2)(1,3)(1,3.5)
\qbezier(1,2)(1,1)(1,.5)
\qbezier(1,2)(.5,2)(0,2)

\qbezier(2.7,2.7)(3,3)(3.5,3.5)
\qbezier(2.7,2.7)(1.3,2.2)(1.3,.7)
\qbezier(2.7,2.7)(3,2.4)(3.5,2.2)

\qbezier(2.7,1.3)(3,1)(3.5,.6)
\qbezier(2.7,1.3)(3,1.6)(3.5,1.8)
\qbezier(2.7,1.3)(2.4,1)(2,.6)



\end{picture}
\begin{center}
$Fig.~\ref{terzabif}.9:~The~bifurcation~from~\alpha~to~\beta~and~from~\beta~
to~\gamma~respectively:$
$3^{rd}~possibility$
\end{center}

The intersection is allowed only in cases 1) and 3).
The phase portrait of $\nabla\tilde{f}_x$ for $x\in\delta$ is
equivalent to the one of $\nabla\tilde{f}_x$ for $x\in\beta$.
%
\end{proof}

Consider $(b)$ (see figure \ref{casob}.10).\\
\setlength{\unitlength}{1cm}
\begin{picture}(12,4.5)
\label{casob}

\thinlines
\qbezier(5,2)(2,2)(1,4)
\qbezier(5,2)(2,2)(1,0)
\qbezier(1,4)(2,2)(1,0)

\thicklines
\qbezier(2,0.5)(2.5,2)(2.5,2.5)
\qbezier(0,3.5)(.3,3)(1.7,3)

\thinlines
\qbezier(11,2)(8,2)(7,4)
\qbezier(11,2)(8,2)(7,0)
\qbezier(7,4)(8,2)(7,0)

\thicklines
\qbezier(8,0.5)(7.7,2)(7.7,3)
\qbezier(6,3.5)(6.3,3)(8.5,2.5)

\put(1.8,1.9){$\alpha$}
\put(2.9,1.9){$\beta$}
\put(1.4,3.1){$\gamma$}

\put(7.5,1.9){$\alpha$}
\put(8.7,1.8){$\beta$}
\put(7.4,3.0){$\gamma$}
\put(7.9,2.7){$\delta$}


\end{picture}
\begin{center}
$Fig.~\ref{casob}.10:~Intersection~of~bifurcation~lines:~case~(b)$
\end{center}

\begin{proposition}
\label{propcasob}
The intersection of bifurcation lines of $(b)$ gives rise to an allowed
bifurcation diagram.
\end{proposition}
\begin{proof}
The phase portrait of $\nabla\tilde{f}_x$, for $x$ in $\alpha$,
exhibits all the gradient lines $\gamma_{ns_i}$, for $i=1,2,3$,
because $\alpha$ is bounded by all the sides $l_i$ of the caustic
(see figure \ref{casob1}.11).\\
\setlength{\unitlength}{1cm}
\begin{picture}(14,5)
\label{casob1}

\thinlines
\put(1,2){\circle*{.1}}
\put(2.7,2.7){\circle*{.1}}
\put(2.7,1.3){\circle*{.1}}
\put(2,2){\circle*{.1}}

\scriptsize
\put(1.1,1.8){$\displaystyle s_1$}
\put(2.2,1.9){$\displaystyle n$}
\put(2.9,2.6){$\displaystyle s_2$}
\put(2.9,1.2){$\displaystyle s_3$}
\normalsize

\qbezier(2,2)(1.5,2)(1,2)
\qbezier(2,2)(2.5,2.5)(2.7,2.7)
\qbezier(2,2)(2.5,1.5)(2.7,1.3)

\qbezier(1,2)(1,3)(1,3.5)
\qbezier(1,2)(1,1)(1,.5)
\qbezier(1,2)(.5,2)(0,2)

\qbezier(2.7,2.7)(3,3)(3.5,3.5)
\qbezier(2.7,2.7)(2.4,3)(2,3.4)
\qbezier(2.7,2.7)(3,2.4)(3.5,2.2)

\qbezier(2.7,1.3)(3,1)(3.5,.6)
\qbezier(2.7,1.3)(3,1.6)(3.5,1.8)
\qbezier(2.7,1.3)(2.4,1)(2,.6)

\put(0,.5){$(\alpha)$}

\end{picture}
\begin{center}
$Fig.~\ref{casob1}.11:~The~phase~portrait~of~\nabla\tilde{f}_x~for~x\in\alpha$
\end{center}
The phase portrait of $\nabla\tilde{f}_x$, for $x$ in
$\beta$ and for $x$ in $\gamma$, are represented in
figure below: from $\alpha$ to $\beta$ the gradient line
$\gamma_{s_2s_1}$ appears and $\gamma_{ns_1}$ breaks; 
from $\alpha$ to $\gamma$ the gradient line
$\gamma_{s_2s_3}$ appears and $\gamma_{ns_3}$ breaks (see figure 
\ref{casob2}.12).\\
\setlength{\unitlength}{1cm}
\begin{picture}(14,5)
\label{casob2}

\thinlines
\put(1,2){\circle*{.1}}
\put(2.7,2.7){\circle*{.1}}
\put(2.7,1.3){\circle*{.1}}
\put(2,2){\circle*{.1}}

\scriptsize
\put(1.1,1.8){$\displaystyle s_1$}
\put(2.2,1.9){$\displaystyle n$}
\put(2.9,2.6){$\displaystyle s_2$}
\put(2.9,1.2){$\displaystyle s_3$}
\normalsize

\qbezier(2,3.2)(1.5,2.5)(1,2)
\qbezier(2,2)(2.5,2.5)(2.7,2.7)
\qbezier(2,2)(2.5,1.5)(2.7,1.3)

\qbezier(1,2)(1,3)(1,3.5)
\qbezier(1,2)(1,1)(1,.5)
\qbezier(1,2)(.5,2)(0,2)

\qbezier(2.7,2.7)(3,3)(3.5,3.5)
\qbezier(2.7,2.7)(1.3,2.2)(1.3,.7)
\qbezier(2.7,2.7)(3,2.4)(3.5,2.2)

\qbezier(2.7,1.3)(3,1)(3.5,.6)
\qbezier(2.7,1.3)(3,1.6)(3.5,1.8)
\qbezier(2.7,1.3)(2.4,1)(2,.6)

\thinlines
\put(7,2){\circle*{.1}}
\put(8.7,2.7){\circle*{.1}}
\put(8.7,1.3){\circle*{.1}}
\put(8,2){\circle*{.1}}

\scriptsize
\put(7.1,1.8){$\displaystyle s_1$}
\put(8.2,1.9){$\displaystyle n$}
\put(8.9,2.7){$\displaystyle s_2$}
\put(8.9,1.2){$\displaystyle s_3$}
\normalsize

\qbezier(8,2)(7.5,2)(7,2)
\qbezier(8,2)(8.5,2.5)(8.7,2.7)

\qbezier(7,2)(7,3)(7,3.5)
\qbezier(7,2)(7,1)(7,.5)
\qbezier(7,2)(6.5,2)(6,2)

\qbezier(8.7,2.7)(9,3)(9.5,3.5)
\qbezier(8.7,2.7)(8.4,3)(8,3.4)
\qbezier(8.7,2.7)(8.4,1.2)(7.7,.8)

\qbezier(8.7,1.3)(9,1)(9.5,.6)
\qbezier(8.7,1.3)(9,1.6)(9.5,1.8)
\qbezier(8.7,1.3)(9,1.9)(9.5,2.3)
\qbezier(8.7,1.3)(8.4,1)(8,.6)

\put(0,.5){$(\beta)$}
\put(6,.5){$(\gamma)$}


\end{picture}\\
\begin{center}
$Fig.~\ref{casob2}.12:~The~phase~portrait~of~\nabla\tilde{f}_x~for~x\in\beta~
and~x\in\gamma~respectively$
\end{center}
The bifurcations
from $\beta$ to $\delta$ and from $\gamma$ to $\delta$, in principle,
can occur into
two different ways, as explained in section (4.5) of \cite{M}.
However, the way gradient lines
$\gamma_{ns_3}$ and $\gamma_{ns_1}$ wind around the node to provide
the bifurcations is fixed by the phase portrait of $\nabla\tilde{f}_t$,
where $t$ is the intersection point of the bifurcation lines. This is
in turn determined by the bifurcations from $\alpha$ to $\beta$
and from $\alpha$ to $\gamma$, and showed in figure \ref{2casob1}.13.\\
%
%
%
%
%
%
%
%
%
%
%
%
%
%
%
%
\setlength{\unitlength}{1cm}
\begin{picture}(14,5)
\label{2casob1}

\thinlines
\put(1,2){\circle*{.1}}
\put(2.7,2.7){\circle*{.1}}
\put(2.7,1.3){\circle*{.1}}
\put(2,2){\circle*{.1}}

\scriptsize
\put(1.1,1.8){$\displaystyle s_1$}
\put(2.2,1.9){$\displaystyle n$}
\put(2.9,2.6){$\displaystyle s_2$}
\put(2.9,1.2){$\displaystyle s_3$}
\normalsize

\qbezier(2.7,2.7)(1,2)(1,2)
\qbezier(2,2)(2.5,2.5)(2.7,2.7)
\qbezier(2.7,2.7)(2.7,1.5)(2.7,1.3)

\qbezier(1,2)(1,3)(1,3.5)
\qbezier(1,2)(1,1)(1,.5)
\qbezier(1,2)(.5,2)(0,2)

\qbezier(2.7,2.7)(3,3)(3.5,3.5)

\qbezier(2.7,1.3)(3,1)(3.5,.6)
\qbezier(2.7,1.3)(3,1.6)(3.5,1.8)
\qbezier(2.7,1.3)(2.4,1)(2,.6)


\end{picture}
\begin{center}
$Fig.~\ref{2casob1}.13:~The~phase~portrait~of~\nabla\tilde{f}_t~where~t~is
~the~intersection~point$
$of~bifurcation~lines$
\end{center}
The phase portrait of $\nabla\tilde{f}_x$, for $x$ in
$\delta$, is unambiguously determined (and shown in figure
\ref{casob3}.14). The intersection is permitted.\\
\setlength{\unitlength}{1cm}
\begin{picture}(14,4)
\label{casob3}

\thinlines
\put(1,2){\circle*{.1}}
\put(2.7,2.7){\circle*{.1}}
\put(2.7,1.3){\circle*{.1}}
\put(2,2){\circle*{.1}}

\scriptsize
\put(1.1,1.8){$\displaystyle s_1$}
\put(2.2,1.9){$\displaystyle n$}
\put(2.9,2.6){$\displaystyle s_2$}
\put(2.9,1.2){$\displaystyle s_3$}
\normalsize

\qbezier(2,3.2)(1.5,2.5)(1,2)
\qbezier(2.7,2.7)(1.3,2.2)(1.3,.7)

\qbezier(1,2)(1,3)(1,3.5)
\qbezier(1,2)(1,1)(1,.5)
\qbezier(1,2)(.5,2)(0,2)

\qbezier(2.7,2.7)(3,3)(3.5,3.5)
\qbezier(2.7,2.7)(2.4,1.2)(1.7,.8)

\qbezier(2.7,1.3)(3,1)(3.5,.6)
\qbezier(2.7,1.3)(3,1.6)(3.5,1.8)
\qbezier(2.7,1.3)(2.4,1)(2,.6)
\qbezier(2.7,1.3)(3,1.9)(3.5,2.3)

\put(0,.5){$(\delta)$}

\end{picture}
\begin{center}
$Fig.~\ref{casob3}.14:~The~phase~portrait~of~\nabla\tilde{f}_x~for~x\in\delta$
\end{center}
%
\end{proof}

Consider now $(c)$ (see figure \ref{casoc1}.15).\\
\setlength{\unitlength}{1cm}
\begin{picture}(12,5.2)
\label{casoc1}

\thinlines
\qbezier(5,2)(2,2)(1,4)
\qbezier(5,2)(2,2)(1,0)
\qbezier(1,4)(2,2)(1,0)

\thicklines
\qbezier(2,0.5)(2.5,2)(2.5,2.5)
\qbezier(2,5)(2,3.5)(1.45,2.8)

\thinlines
\qbezier(11,2)(8,2)(7,4)
\qbezier(11,2)(8,2)(7,0)
\qbezier(7,4)(8,2)(7,0)

\thicklines
\qbezier(8,0.5)(8.5,2)(8.25,2.6)
\qbezier(10,5)(10,2)(7.5,2)

\put(1.8,1.9){$\alpha$}
\put(2.9,1.9){$\beta$}
\put(1.3,3.1){$\gamma$}

\put(7.8,1.6){$\alpha$}
\put(8.7,1.8){$\beta$}
\put(7.7,2.4){$\gamma$}
\put(8.5,2.3){$\delta$}


\end{picture}
\begin{center}
$Fig.~\ref{casoc1}.15:~Intersection~of~bifurcation~lines:~case~(c)$
\end{center}

\begin{proposition}
\label{propcasoc}
Bifurcation lines of case $(c)$ can intersect, however the new
bifurcation diagram is allowed provided it also contains
a new bifurcation line, from the intersection point to $l_2$ and lying 
inside $\beta$, at whose points the phase portrait of $\nabla\tilde{f}_x$
exhibits a non-generic gradient line $\gamma_{s_1s_3}$ from $s_1$ to $s_3$.
The
bifurcation locus after the intersection is shown in figure \ref{casoc4}.16.\\
\end{proposition}

\setlength{\unitlength}{1cm}
\begin{picture}(14,5)
\label{casoc4}
\thinlines
\qbezier(5,2)(2,2)(1,4)
\qbezier(5,2)(2,2)(1,0)
\qbezier(1,4)(2,2)(1,0)

\thicklines
\qbezier(2,0.5)(2.5,2)(2.25,2.6)
\qbezier(4,5)(4,2)(1.5,2)
\qbezier(2.35,2.1)(2.8,2)(3.5,2.15)
\put(1.8,1.6){$\alpha$}
\put(2.7,1.7){$\beta$}
\put(1.7,2.4){$\gamma$}
\put(2.5,2.3){$\delta$}
\put(2.9,2.1){$\epsilon$}


\end{picture}
\begin{center}
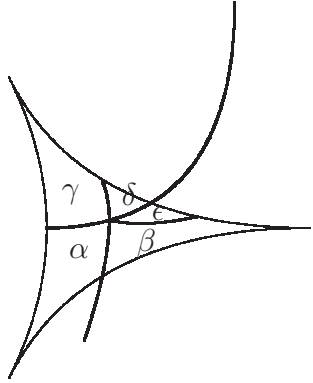
$Fig.~\ref{casoc4}.16:~The~resulting~bifurcation~diagram~after~the~
intersection$
$of~bifurcation~lines~in~case~(c)$
\end{center}

\begin{proof}
Whatever the position of the third bifurcation half-line is, the phase
portrait of $\nabla\tilde{f}_x$, for $x$ in $\alpha$, exhibits all the
gradient lines $\gamma_{ns_i}$, for $i=1,2,3$ (see figure \ref{casoc2}.17).\\
\setlength{\unitlength}{1cm}
\begin{picture}(14,5)
\label{casoc2}

\thinlines
\put(1,2){\circle*{.1}}
\put(2.7,2.7){\circle*{.1}}
\put(2.7,1.3){\circle*{.1}}
\put(2,2){\circle*{.1}}

\scriptsize
\put(1.1,1.8){$\displaystyle s_1$}
\put(2.2,1.9){$\displaystyle n$}
\put(2.9,2.6){$\displaystyle s_2$}
\put(2.9,1.2){$\displaystyle s_3$}
\normalsize

\qbezier(2,2)(1.5,2)(1,2)
\qbezier(2,2)(2.5,2.5)(2.7,2.7)
\qbezier(2,2)(2.5,1.5)(2.7,1.3)

\qbezier(1,2)(1,3)(1,3.5)
\qbezier(1,2)(1,1)(1,.5)
\qbezier(1,2)(.5,2)(0,2)

\qbezier(2.7,2.7)(3,3)(3.5,3.5)
\qbezier(2.7,2.7)(2.4,3)(2,3.4)
\qbezier(2.7,2.7)(3,2.4)(3.5,2.2)

\qbezier(2.7,1.3)(3,1)(3.5,.6)
\qbezier(2.7,1.3)(3,1.6)(3.5,1.8)
\qbezier(2.7,1.3)(2.4,1)(2,.6)

\put(0,.5){$(\alpha)$}

\end{picture}
\begin{center}
$Fig.~\ref{casoc2}.17:~The~phase~portrait~of~\nabla\tilde{f}_x~for~x\in\alpha$
\end{center}
There is only one way
to perform the bifurcations from $\alpha$ to $\beta$ and from $\alpha$
to $\gamma$, obtaining the phase portrait showed in figure \ref{casoc3}.18.\\
\setlength{\unitlength}{1cm}
\begin{picture}(14,5)
\label{casoc3}

\thinlines
\put(1,2){\circle*{.1}}
\put(2.7,2.7){\circle*{.1}}
\put(2.7,1.3){\circle*{.1}}
\put(2,2){\circle*{.1}}

\scriptsize
\put(1.1,1.8){$\displaystyle s_1$}
\put(2.2,1.9){$\displaystyle n$}
\put(2.9,2.6){$\displaystyle s_2$}
\put(2.9,1.2){$\displaystyle s_3$}
\normalsize

\qbezier(2,3.2)(1.5,2.5)(1,2)
\qbezier(2,2)(2.5,2.5)(2.7,2.7)
\qbezier(2,2)(2.5,1.5)(2.7,1.3)

\qbezier(1,2)(1,3)(1,3.5)
\qbezier(1,2)(1,1)(1,.5)
\qbezier(1,2)(.5,2)(0,2)

\qbezier(2.7,2.7)(3,3)(3.5,3.5)
\qbezier(2.7,2.7)(1.3,2.2)(1.3,.7)
\qbezier(2.7,2.7)(3,2.4)(3.5,2.2)

\qbezier(2.7,1.3)(3,1)(3.5,.6)
\qbezier(2.7,1.3)(3,1.6)(3.5,1.8)
\qbezier(2.7,1.3)(2.4,1)(2,.6)

\thinlines
\put(7,2){\circle*{.1}}
\put(8.7,2.7){\circle*{.1}}
\put(8.7,1.3){\circle*{.1}}
\put(8,2){\circle*{.1}}

\scriptsize
\put(7.1,1.8){$\displaystyle s_1$}
\put(8.2,1.9){$\displaystyle n$}
\put(8.9,2.7){$\displaystyle s_2$}
\put(8.9,1.2){$\displaystyle s_3$}
\normalsize

\qbezier(8,1.5)(7,1.5)(7,2)
\qbezier(9.5,2.7)(8.7,1.5)(8,1.5)
\qbezier(8,2)(8.5,2.5)(8.7,2.7)

\qbezier(7,2)(7,3)(7,3.5)
\qbezier(7,2)(7.5,2)(8,2)
\qbezier(7,2)(6.5,2)(6,2)

\qbezier(8.7,2.7)(9,3)(9.5,3.5)
\qbezier(8.7,2.7)(8.4,3)(8,3.4)

\qbezier(8.7,1.3)(9,1)(9.5,.6)
\qbezier(8.7,1.3)(9,1.6)(9.5,1.8)
\qbezier(8.7,1.3)(8.4,1)(7,.6)
\qbezier(8.7,1.3)(9,1.9)(9.5,2.3)

\put(0,0.5){$(\beta)$}
\put(6,0.5){$(\gamma)$}


\end{picture}
\begin{center}
$Fig.~\ref{casoc3}.18:~The~phase~portrait~of~\nabla\tilde{f}_x~for~x\in\beta~
and~x\in\gamma~respectively$
\end{center}
According to proposition (4.19) in \cite{M},
in the phase portrait of $\nabla\tilde{f}_x$ for $x\in\beta$, shown 
in figure \ref{casoc3}.18, a gradient
line from $s_1$ to $s_3$ would provide a non-allowed diagram, so 
$\beta$ can not be bounded
by the bifurcation line, shown in figure \ref{casoc1}.15, separating 
$\beta$ from
$\delta$.
In $\gamma$, instead, the gradient line $\gamma_{s_2s_1}$ gives an allowed
bifurcation diagram,
so 
the bifurcation line separating $\gamma$ from 
$\delta$ is allowed. 
The
bifurcation line separating $\delta$ from 
$\beta$ produces an allowed bifurcation
diagram, 
because the gradient line
$\gamma_{s_1s_3}$ can appear in phase portrait $\delta$, 
though this implies for $x\in\beta$ a phase diagram of $\nabla\tilde{f}_x$
different from that shown in figure \ref{casoc3}.18.
However, as already explained in subsection (4.6) of \cite{M},
we can suppose that the
intersection point $t\in{\cal B}_{(2,1),(1,3)}$ belongs also to 
${\cal B}_{2,3}$, so a
new bifurcation line ${\cal B'}$ arises: ${\cal B'}$ is actually a segment 
with an
extreme at the intersection point $t$ and  
the other extreme on $l_2$,
where the saddle $s_2$ glues together with the node.
At a point $x\in{\cal B'}$, a saddle-to-saddle separatrix 
$\gamma_{s_2s_3}$ appears in the phase portrait of
$\nabla\tilde{f}_x$. 
Since the problem is crossing the bifurcartion line
from $\beta$ to $\delta$, we suppose that ${\cal B'}$
lies inside $\beta$.
Call $\epsilon$ the new subset determined in the bifurcation
diagram by ${\cal B'}$, lying between $\beta$ and $\delta$. 
At points of ${\cal B'}$, from $\beta$ to $\epsilon$, 
in principle there are two ways to
perform the bifurcation, that is,
$\gamma_{s_2s_3}$ can be obtained by joining two different pairs of
separatrices, however the choice is fixed, and shown in figure 
\ref{casoc5}.19, by the phase
portrait of $\nabla\tilde{f}_t$, where $t$ is the intersection point of
bifurcation lines (see figure \ref{casoc8}.22).\\
\setlength{\unitlength}{1cm}
\begin{picture}(14,4.5)
\label{casoc5}

\thinlines
\put(1,2){\circle*{.1}}
\put(2.7,2.7){\circle*{.1}}
\put(2.7,1.3){\circle*{.1}}
\put(2,2){\circle*{.1}}

\scriptsize
\put(1.1,1.8){$\displaystyle s_1$}
\put(2.2,1.9){$\displaystyle n$}
\put(2.9,2.6){$\displaystyle s_2$}
\put(2.9,1.2){$\displaystyle s_3$}
\normalsize

\qbezier(2,3.2)(1.5,2.5)(1,2)
\qbezier(2,2)(2.5,2.5)(2.7,2.7)

\qbezier(1,2)(1,3)(1,3.5)
\qbezier(1,2)(1,1)(1,.5)
\qbezier(1,2)(.5,2)(0,2)

\qbezier(2.7,2.7)(3,3)(3.5,3.5)
\qbezier(2.7,2.7)(1.5,2.7)(1.5,2)
\qbezier(1.5,2)(1.5,1.3)(2.7,1.3)
\qbezier(2.7,2.7)(3,2.4)(3.5,2.2)

\qbezier(2.7,1.3)(3,1)(3.5,.6)
\qbezier(2.7,1.3)(3,1.6)(3.5,1.8)
\qbezier(2.7,1.3)(2.4,1)(2,.6)

\thinlines
\put(7,2){\circle*{.1}}
\put(8.7,2.7){\circle*{.1}}
\put(8.7,1.3){\circle*{.1}}
\put(8,2){\circle*{.1}}

\scriptsize
\put(7.1,1.8){$\displaystyle s_1$}
\put(8.2,1.9){$\displaystyle n$}
\put(8.95,2.75){$\displaystyle s_2$}
\put(8.9,1.15){$\displaystyle s_3$}
\normalsize

\qbezier(8,3.5)(7.5,2.5)(7,2)
\qbezier(8,2)(8.5,2.5)(8.7,2.7)

\qbezier(7,2)(7,3)(7,3.7)
\qbezier(7,2)(7,.5)(8,.5)
\qbezier(7,2)(6.5,2)(6,2)

\qbezier(8.7,2.7)(9,3)(9.5,3.5)
\qbezier(8.7,2.7)(7.7,2.7)(7.7,2)
\qbezier(7.7,2)(7.7,1.3)(8.7,1.5)
\qbezier(8.7,1.5)(9.5,1.7)(9.5,1.9)
\qbezier(8.7,2.7)(9,2.7)(9.5,2.7)

\qbezier(8.7,1.3)(9,1)(9.5,.6)
\qbezier(8.7,1.3)(9,1.3)(9.5,1.7)
\qbezier(8.7,1.3)(7.4,.5)(7.4,2)
\qbezier(7.4,2)(7.7,2.8)(8.5,3.2)
\qbezier(8.7,1.3)(8.9,.8)(9.2,.6)

\put(0,0){$(\beta\rightarrow\epsilon)$}
\put(6,0){$(\epsilon)$}


\end{picture}
\begin{center}
$Fig.~\ref{casoc5}.19:~The~bifurcation~from~\beta~to~\epsilon~and~
the~phase~portrait~of~\nabla\tilde{f}_x$
$for~x\in\epsilon$
\end{center}
The bifurcation from $\epsilon$ to
$\delta$, being characterized by a non-generic gradient line
$\gamma_{s_1s_3}$ in the phase portrait of $\nabla\tilde{f}_x$, 
is allowed. Since such line
appears in the phase portrait of $\nabla\tilde{f}_t$, where
$t$ is the intersection point of bifurcation lines, the choice of the
pair of separatices to be joined in $\gamma_{s_1s_3}$ is 
thus determined (see figure
\ref{casoc6}.20).\\
\setlength{\unitlength}{1cm}
\begin{picture}(14,5)
\label{casoc6}

\thinlines
\put(1,2){\circle*{.1}}
\put(2.7,2.7){\circle*{.1}}
\put(2.7,1.3){\circle*{.1}}
\put(2,2){\circle*{.1}}

\scriptsize
\put(1.1,1.8){$\displaystyle s_1$}
\put(2.2,1.9){$\displaystyle n$}
\put(2.95,2.75){$\displaystyle s_2$}
\put(2.9,1.15){$\displaystyle s_3$}
\normalsize

\qbezier(2,3.5)(1.5,2.5)(1,2)
\qbezier(2,2)(2.5,2.5)(2.7,2.7)

\qbezier(1,2)(1,3)(1,3.7)
\qbezier(1,2)(1,1.3)(2.7,1.3)
\qbezier(1,2)(.5,2)(0,2)

\qbezier(2.7,2.7)(3,3)(3.5,3.5)
\qbezier(2.7,2.7)(1.7,2.7)(1.7,2)
\qbezier(1.7,2)(1.7,1.3)(2.7,1.5)
\qbezier(2.7,1.5)(3.5,1.7)(3.5,1.9)
\qbezier(2.7,2.7)(3,2.7)(3.5,2.7)

\qbezier(2.7,1.3)(3,1)(3.5,.6)
\qbezier(2.7,1.3)(3,1.3)(3.5,1.7)
\qbezier(2.7,1.3)(2.7,.6)(2.2,.6)

\thinlines
\put(7,2){\circle*{.1}}
\put(8.7,2.7){\circle*{.1}}
\put(8.7,1.3){\circle*{.1}}
\put(8,2){\circle*{.1}}

\scriptsize
\put(7,1.7){$\displaystyle s_1$}
\put(8.2,1.9){$\displaystyle n$}
\put(8.95,2.75){$\displaystyle s_2$}
\put(8.9,1.15){$\displaystyle s_3$}
\normalsize

\qbezier(8,3.5)(7.5,2.5)(7,2)
\qbezier(8,2)(8.5,2.5)(8.7,2.7)

\qbezier(7,2)(7,3)(7,3.7)
\qbezier(7,2)(8.7,1)(9.5,1.9)
\qbezier(7,2)(6.5,2)(6,2)

\qbezier(8.7,2.7)(9,3)(9.5,3.5)
\qbezier(8.7,2.7)(7.7,2.7)(7.7,2)
\qbezier(7.7,2)(7.7,1.7)(8.7,1.8)
\qbezier(8.7,1.8)(9.5,1.9)(9.5,2.2)
\qbezier(8.7,2.7)(9,2.7)(9.5,2.7)

\qbezier(8.7,1.3)(9,1)(9.5,.6)
\qbezier(8.7,1.3)(9,1.3)(9.5,1.7)
\qbezier(8.7,1.3)(8.7,.6)(8.2,.6)
\qbezier(8.7,1.3)(7,1.4)(6,1.5)

\put(0,0){$(\epsilon\rightarrow\delta)$}
\put(6,0){$(\delta)$}


\end{picture}
\begin{center}
$Fig.~\ref{casoc5}.20:~The~bifurcation~from~\epsilon~to~\delta~and~
the~phase~portrait~of~\nabla\tilde{f}_x$
$for~x\in\delta$
\end{center}
For $x$ belonging to the bifurcation line separating $\delta$ from 
$\gamma$, a non-generic gradient line $\gamma_{s_2s_1}$
must appear in the phase portrait of $\nabla\tilde{f}_x$: 
again, since such a line
appears in the phase portrait of $\nabla\tilde{f}_t$,
the choice of the pair of separatrices intersecting in
$\gamma_{s_2s_1}$
is fixed (see figure
\ref{casoc7}.21).\\  
\setlength{\unitlength}{1cm}
\begin{picture}(14,5)
\label{casoc7}

\thinlines
\put(1,2){\circle*{.1}}
\put(2.7,2.7){\circle*{.1}}
\put(2.7,1.3){\circle*{.1}}
\put(2,2){\circle*{.1}}

\scriptsize
\put(1,1.7){$\displaystyle s_1$}
\put(2.2,1.9){$\displaystyle n$}
\put(2.9,2.75){$\displaystyle s_2$}
\put(2.9,1.15){$\displaystyle s_3$}
\normalsize

\qbezier(2.7,2.7)(1.85,2.35)(1,2)
\qbezier(2,2)(2.5,2.5)(2.7,2.7)

\qbezier(1,2)(1,3)(1,3.7)
\qbezier(1,2)(2.7,1)(3.5,1.9)
\qbezier(1,2)(.5,2)(0,2)

\qbezier(2.7,2.7)(3,3)(3.5,3.5)
\qbezier(2.7,2.7)(3,2.7)(3.5,2.7)

\qbezier(2.7,1.3)(3,1)(3.5,.6)
\qbezier(2.7,1.3)(3,1.3)(3.5,1.7)
\qbezier(2.7,1.3)(2.7,.6)(2.2,.6)
\qbezier(2.7,1.3)(1,1.4)(0,1.5)

\thinlines
\put(7,2){\circle*{.1}}
\put(8.7,2.7){\circle*{.1}}
\put(8.7,1.3){\circle*{.1}}
\put(8,2){\circle*{.1}}

\scriptsize
\put(7,1.7){$\displaystyle s_1$}
\put(8.2,1.9){$\displaystyle n$}
\put(8.9,2.75){$\displaystyle s_2$}
\put(8.9,1.15){$\displaystyle s_3$}
\normalsize

\qbezier(8,2)(7.5,2)(7,2)
\qbezier(8,2)(8.5,2.5)(8.7,2.7)

\qbezier(7,2)(7,3)(7,3.7)
\qbezier(7,2)(8.7,1)(9.5,1.9)
\qbezier(7,2)(6.5,2)(6,2)

\qbezier(8.7,2.7)(9,3)(9.5,3.5)
\qbezier(8.7,2.7)(9,2.7)(9.5,2.7)
\qbezier(8.7,2.7)(8.4,3)(8,3.4)

\qbezier(8.7,1.3)(9,1)(9.5,.6)
\qbezier(8.7,1.3)(9,1.3)(9.5,1.7)
\qbezier(8.7,1.3)(8.7,.6)(8.2,.6)
\qbezier(8.7,1.3)(7,1.4)(6,1.5)

\put(0,0){$(\delta\rightarrow\gamma)$}
\put(6,0){$(\gamma)$}


\end{picture}
\begin{center}
$Fig.~\ref{casoc5}.21:~The~bifurcation~from~\delta~to~\gamma~and~
the~phase~portrait~of~\nabla\tilde{f}_x$
$for~x\in\gamma$
\end{center}
The phase portrait of $\nabla\tilde{f}_t$, where
$t$ is the intersection point, is shown in figure \ref{casoc8}.22.\\
\setlength{\unitlength}{1cm}
\begin{picture}(14,4.5)
\label{casoc8}

\thinlines
\put(1,2){\circle*{.1}}
\put(2.7,2.7){\circle*{.1}}
\put(2.7,1.3){\circle*{.1}}
\put(2,2){\circle*{.1}}

\scriptsize
\put(.7,1.75){$\displaystyle s_1$}
\put(2.2,1.9){$\displaystyle n$}
\put(2.9,2.6){$\displaystyle s_2$}
\put(2.9,1.2){$\displaystyle s_3$}
\normalsize

\qbezier(2.7,2.7)(1.85,2.35)(1,2)
\qbezier(2,2)(2.5,2.5)(2.7,2.7)
\qbezier(1,2)(1.85,1.65)(2.7,1.3)

\qbezier(1,2)(1,3)(1,3.5)

\qbezier(1,2)(.5,2)(0,2)

\qbezier(2.7,2.7)(3,3)(3.5,3.5)

\qbezier(2.7,2.7)(3,2.4)(3.5,2.2)

\qbezier(2.7,1.3)(3,1)(3.5,.6)
\qbezier(2.7,1.3)(3,1.6)(3.5,1.8)
\qbezier(2.7,1.3)(2.4,1)(2,.6)


\end{picture}
\begin{center}
$Fig.~\ref{casoc8}.22:~The~phase~portrait~of~\nabla\tilde{f}_t~where~t
~is~the~intersection~point$
$of~bifurcation~lines$
\end{center}
%
\end{proof}

\begin{remark}
\rm
Observe that there are five ways to break the exceptional gradient
lines appearing in the phase portrait of $\nabla\tilde{f}_t$, shown in
figure \ref{casoc8}.22, as five are the bifurcation lines arising
from the intersection point in the bifurcation diagram in figure
\ref{casoc1}.16. Note that $\gamma_{s_2s_1}$
and $\gamma_{s_1s_3}$, as already explained in subsection (4.6) of \cite{M},
can break in such a way to give rise to the exceptional
gradient line $\gamma_{s_2s_3}$, which is not exhibited by any of phase
portraits of $\nabla\tilde{f}_x$, for $x$ belonging to the bifurcation
lines, when these do not intersect.
\end{remark}

Consider now $(d)$ (see figure \ref{casod}.23). 
Among those in figure \ref{allowedandnot}.5, 
the only allowed diagram exhibiting this
configuration is (D).\\
\setlength{\unitlength}{1cm}
\begin{picture}(12,5.5)
\label{casod}
\thinlines
\qbezier(5,2)(2,2)(1,4)
\qbezier(5,2)(2,2)(1,0)
\qbezier(1,4)(2,2)(1,0)

\thicklines
\qbezier(2,0.5)(2.5,2)(2.5,2.5)
\qbezier(2,5)(3,5)(3,1.7)


\end{picture}
\begin{center}
$Fig.~\ref{casod}.23:~Intersection~of~bifurcation~lines:~case~(d)$
\end{center}

\begin{proposition}
\label{propcasod}
The bifurcation diagram resulting from the intersection of 
bifurcation lines of $(d)$ is not allowed.
\end{proposition}

\begin{proof}
If the intersection were allowed, then the phase portrait of 
$\nabla\tilde{f}_t$, where $t$ is the
intersection point, would exhibit 
the gradient lines $\gamma_{ns_2}$, $\gamma_{ns_3}$
and the non-generic gradient lines $\gamma_{s_2s_1}$ and
$\gamma_{s_3s_1}$, but such a vector field does not exist.
\end{proof}

Consider $(e)$ (see figure \ref{casoe}.24).\\
\setlength{\unitlength}{1cm}
\begin{picture}(12,4.5)
\label{casoe}

\thinlines
\qbezier(5,2)(2,2)(1,4)
\qbezier(5,2)(2,2)(1,0)
\qbezier(1,4)(2,2)(1,0)

\thicklines
\qbezier(4.5,1)(3.5,1)(3.5,2.15)
\qbezier(2,0.2)(2,.7)(1.3,0.7)


\end{picture}
\begin{center}
$Fig.~\ref{casod}.24:~Intersection~of~bifurcation~lines:~case~(e)$
\end{center}

\begin{proposition}
\label{propcasoe}
The bifurcation diagram resulting from the intersection of the
bifurcation lines of $(e)$ is not allowed.
\end{proposition}
\begin{proof}
If the intersection were allowed, 
the phase portrait of $\nabla\tilde{f}_t$, when $t$ is the
intersection point, would exhibit two non-generic gradient lines
$\gamma_{s_1s_2}$ and $\gamma_{s_2s_1}$, giving a contraddiction. 
\end{proof}

We realize that, inside the caustic, 
the bifurcation locus of a perturbation of the
elliptic umbilic can be rather complex.
We can resume all the results in the following theorem:

\begin{theorem}
\label{finalresult}
The bifurcation locus ${\cal B}$ of a small perturbation of the 
generating function
(\ref{ellipticumbilic}) of the elliptic umbilic 
in dimension 2 has the following
features:\\
- outside a compact subset containing 
the caustic $K$, ${\cal B}$ exhibits three bifurcation half-lines;\\
- generically, these half-lines intersect $K$ along one of
its sides and have their 
extreme at a fold point of one of the opposite sides;\\
- the allowed mutual positions of bifurcation lines and their
possible intersections are described by lemma \ref{lemdiagr} and
by propositions \ref{propcasoa}, \ref{propcasob}, \ref{propcasoc},
\ref{propcasod}, \ref{propcasoe};\\
- inside $K$, ${\cal B}$ can also contain 
immersed $S^1$'s or segments, with extremes on the same side of $K$.
\end{theorem}


\begin{thebibliography}{100}\frenchspacing\small

%
%
%
%
%
%
%
%
\bibitem{F2} K. Fukaya, \emph{Multivalued Morse theory, asymptotics
analysis and mirror symmetry,} (2002). Available from the web page
{\tt http://www.math.kyoto-u.ac.jp/\~{}fukaya/fukayagrapat.dvi} 

%
%
%
%
\bibitem{M} G. Marelli, \emph{Two-dimensional Lagrangian singularities 
and bifurcations of gradient lines I,} pre-print.
%
%
%
%

\end{thebibliography}
\end{document}